\documentclass[a4paper,11pt]{amsart}

\usepackage[english]{babel}
\usepackage[centertags]{amsmath}
\usepackage{latexsym, amsfonts, amssymb, graphicx}

\usepackage[all,dvips]{xy}
\usepackage{epic}
\usepackage{xypic}
\usepackage[latin1]{inputenc}
\usepackage{url}
\usepackage{enumerate}

\usepackage{epsfig}
\usepackage{psfrag}
\usepackage{graphicx}

\theoremstyle{plain}
\newtheorem{thm}{Theorem}[section]
\newtheorem{lem}[thm]{Lemma}
\newtheorem{prop}[thm]{Proposition}
\newtheorem{cor}[thm]{Corollary}

\theoremstyle{definition}
\newtheorem{defn}[thm]{Definition}
\newtheorem{nota}[thm]{Notation}

\theoremstyle{remark}
\newtheorem{rem}[thm]{Remark}

\newtheorem{exmp}[thm]{Example}

\def\P{{\mathbf P}}

\newcommand{\F}{\mathbf{F}}
\newcommand{\supp}{\textrm{Supp}}
\newcommand{\ev}{\textrm{ev}}
\newcommand{\Ff}{\mathcal{F}}

\begin{document}

\title[The dual minimum distance of arbitrary dimensional AG codes]{The dual minimum distance of arbitrary--dimensional algebraic--geometric codes}

\author{Alain \textsc{Couvreur}}

%\date{\today}

\address{INRIA Saclay - École polytechnique,
Laboratoire d'informatique (LIX),
UMR 7161,
91128 Palaiseau Cedex,
France}
\email{alain.couvreur@inria.fr}

\begin{abstract}
In this article, the minimum distance of the dual $C^{\bot}$ of a functional code $C$ on an arbitrary--dimensional variety $X$ over a finite field $\F_q$ is studied.
The approach is based on problems {\it \`a la Cayley--Bacharach} and consists in describing the minimal configurations of points on $X$ which fail to impose independent conditions on forms of some fixed degree $m$.
If $X$ is a curve, the result improves in some situations the well-known \textit{Goppa designed distance}.
\end{abstract}

\maketitle

\thispagestyle{empty}

\noindent \textbf{AMS Classification:} 14J20, 94B27, 14C20.

\noindent \textbf{Keywords:} Algebraic geometry, finite fields, error--correcting codes, algebraic--geometric codes, linear systems.

\section*{Introduction}

A classical problem in  coding theory is the estimation of the minimum distance of some code or family of codes constructed on some variety or some family of varieties.
For algebraic--geometric codes on curves, one easily gets such a lower bound, frequently called the \textit{Goppa designed distance} (see \cite[Definition II.2.4]{sti}).

On higher--dimensional varieties, the problem becomes really harder even when the geometry of the involved variety is well understood. 
This difficulty can be explained by a citation from Little in the introduction of a survey on the topic \cite[Chapter 7]{bouquindiego}: ``the first major difference between higher--dimensional varieties and curves is that points on $X$ of dimension $\geq 2$ are [...] not divisors''.
Therefore, if getting the Goppa designed minimum distance is an easy exercise of function fields theory, obtaining any relevant information on the minimum distance of an algebraic--geometric code on a higher--dimensional variety (or a family of varieties) is often the purpose of an entire article.
For instance, codes on quadrics are studied in \cite{aubry}, some general bounds on codes on arbitrary--dimensional varieties are given in \cite{lachaud} and, in \cite{zarzar}, codes on surfaces having a low Neron--Severi rank are studied (the list is far from being exhaustive).

Another kind of codes associated to algebraic varieties can be studied: \textit{the dual of a functional code}. That is, the orthogonal space for the canonical inner product in $\F_q^n$.
On a curve $X$, the dual of a functional code is also a functional code on $X$ (see \cite[Proposition II.2.10]{sti}).
It turns out that this result does not hold for higher--dimensional varieties.
Such a difference with codes on curves has been felt by Voloch and Zarzar who noticed it in \cite{agctvoloch} and then proved in \cite[\S 10]{oim} using an elementary example of surface (or a higher--dimensional variety, see \cite[Remark II.5.5]{thesebibi}).

Therefore, on varieties of dimension greater than or equal to $2$, one can say that \textit{a new class of codes appears} and it is natural to wonder if this new class contains good codes.
This motivates the study of the parameters of these duals of functional codes on arbitrary--dimensional varieties, which is the purpose of this article.

% In general, the problem of estimating or lower bounding the minimum distance(s) of some code or family of codes is very difficult.
% For some specific families of codes (for instance Reed--Muller codes or algebraic--geometric codes), this problem can be translated into a geometrical or arithmetical problem.

In the present paper, we translate the problem of finding the dual minimum distance of an algebraic--geometric code into a problem of finding some particular configurations of points in a projective space.
In particular, we introduce the elementary notion of \textit{minimally $m$--linked points} (Definition \ref{mlink}), that is sets of points which fail to impose independent conditions on forms of degree $m$ and are minimal for this property.
This notion relates to problems {\it à la Cayley--Bacharach} (see \cite{EGH}) and is central for the proof of Theorem \ref{main}, which gives estimates or lower bounds for the minimum distance of the duals of functional codes.
From a more geometrical point of view, we give the complete description of minimally $m$--linked configurations of less than $3m$ points in any projective space.
It is stated in \cite{EGH} that complete intersections provide such configurations. In addition, the authors ask whether these configurations are the only ones. We give a positive answer to this question for configurations of cardinality lower than or equal to $3m$.

From the coding theoretic point of view, the most surprising application of this result is the case when the variety is a plane curve.
Indeed, in this situation, since the dual of an algebraic--geometric code on a curve is also an algebraic--geometric code on this curve, the dual minimum distance has a lower bound given by the Goppa designed distance.
Therefore, we compare the bound yielded by Theorem \ref{main} with the Goppa designed distance.
It turns out that our bound is better than Goppa's one in two situations.
First, when Goppa's bound is negative and hence irrelevant, since our bound is always positive.
Second, if one can check some incidence condition on the points of evaluation, in this second situation, one can get a bound which is much better than that of Goppa.

Some proofs of the present paper are long and need the treatment of numerous cases. This is the reason why we chose to study examples of applications of the results (in Section \ref{secex}) before proving them.
The study of configurations of points and linear systems having prescribed points in their base locus is often very technical.
For instance see the proof of \cite[Proposition V.4.3]{H}. 

\subsection*{Contents}
Section \ref{secAG} is a brief review on algebraic--geometric codes on curves and arbitrary--dimensional varieties.
Section \ref{secmg} is devoted to the definition of the notion of \textit{$m$--general} and \textit{minimally $m$--linked} configurations of points in a projective space.
The connection between this notion and the dual minimum distance is explained at the beginning of Section \ref{secdual}. 
In addition, Section \ref{secdual} contains the main theorem (Theorem \ref{main}) and its ``geometric version'' (Theorem \ref{geomain}).
Theorem \ref{main} gives lower bounds for the minimum distance of the dual of a functional code.
Explicit examples of applications of the main theorem are presented in Section \ref{secex}.
In particular the case of codes on plane curves and the improvements of the Goppa designed distance are studied.

Sections \ref{secmlink} to \ref{secfourth} are devoted to the proof of Theorem \ref{geomain}.
In Section \ref{secmlink}, two key tools for this proof, namely Lemma \ref{indmin} and Theorem \ref{intersec} are stated.
Lemma \ref{indmin} is a useful trick to handle minimally $m$--linked configurations of points and Theorem \ref{intersec} is one of the numerous formulations of Cayley--Bacharach theorem.
Afterwards, Sections \ref{secfirst} to \ref{secfourth} are devoted to the proofs of some results on configurations of points in projective spaces, yielding the proof of Theorem \ref{main}.

\section{algebraic--geometric codes}\label{secAG}
Let $X$ be a smooth geometrically connected projective variety defined over a finite field $\F_q$. Let $G$ be a divisor on $X$ and $P_1, \ldots , P_n$ be a family of rational points of $X$ avoiding the support of $G$. Denote by $\Delta$ the $0$-cycle defined by the formal sum $\Delta:=P_1 + \cdots +P_n$.
In \cite{manin}, Vl\u{a}du\c{t} and Manin define the functional code $C_L (X, \Delta , G)$ to be the image of the map 
$$
\textrm{ev}_{\Delta}: \left\{
  \begin{array}{ccc}
    L(G) & \rightarrow & \F_q^n \\
    f & \mapsto & (f(P_1), \ldots , f(P_n)) 
  \end{array} ,
\right.
$$
where $L(G)$ denotes the Riemann--Roch space associated to $G$.
When there is no possible confusion on the involved variety, one can remove the ``$X$'' and denote this code by $C_L (\Delta,G)$.

As said in the Introduction, the aim of this paper is to study the minimum distance of the dual code $C_L (X, \Delta , G)^{\bot}$.

\medbreak

\noindent \textbf{Caution.} A usual abuse of notation in coding theory consists in using the term ``dual'' to denote the orthogonal space $C^{\bot}$ of a subspace $C$ of $\F_q^n$ for the canonical inner product.
This space differs from the genuine dual $C^{\vee}$ of $C$, which is the space of linear forms on $C$.
According to the conventions in coding theory, we allow ourselves such an abuse of language in this paper, even if actual dual spaces will be also involved sometimes.
The exponents $\bot$ and $\vee$ enable to differentiate one ``dual'' from another, avoiding any confusion.

\section{Points in $m$--general position}\label{secmg}

In the present section, the base field $k$ is arbitrary.

\subsection{General position in the literature}
The notion of ``general position'' is classical in algebraic geometry. However, there does not seem to exist any consensual definition.
Roughly speaking, a fixed number $s$ of points on a variety $X$ is said to be in \textit{general position} if they correspond to a point of a Zariski dense subset of the space of configurations of $s$ points of $X$.
The point is that the involved \textit{dense subset} depends on the problem we are working on.

The most usual definition is that $s$ points of the affine space $\mathbf{A}_k^r$ (resp. the projective space $\P_k^r$) are in general position if for all $l\leq r$, no $l+2$ of them lie on an $l$--dimensional linear subspace. 
However, several different definitions exist in the literature.   
For instance, the definition of Hartshorne in \cite[Exercise V.4.15]{H} differs from that of Mumford in \cite[Lecture 20]{mumf}.

\subsection{Definition in the present paper} The definition we will use involves linear independence of evaluation maps on a space of homogeneous forms of fixed degree.

\begin{nota}
We denote by $\Ff _{m,r}(k)$ the space $H^0 (\P_k^r, \mathcal{O}_{\P_k^r} (m))$ of homogeneous forms of degree $m$ in $r+1$ variables. If there is no possible confusion on the base field, we denote this space by $\Ff_{m,r}$.  
\end{nota}

Notice that the evaluation at a point of $X$ (or a point of $\P_k^r$ actually) does not make sense for homogeneous forms. To avoid this problem, one can choose a system of homogeneous coordinates in $\P_k^r$ and use the evaluation maps defined in \cite{lachaud2}.

\begin{defn}[{\cite[\S 3]{lachaud2}}]\label{evmap}
Let $P=(p_0 : \cdots : p_r)$ be a rational point of $\P_k^r$.
Let $i\in \{0, \ldots , r\}$ be the smallest integer such that $p_i \neq 0$.
For a nonnegative integer $m$ we define the evaluation map to be 
$$
\ev_P : \left\{
  \begin{array}{ccc}
    \mathcal{F}_{m,r} & \rightarrow & k \\
    f & \mapsto & \frac{f(p_0, \ldots , p_r)}{p_i^m}
  \end{array}
\right. .
$$
\end{defn}

\begin{rem}
The previous definition can be regarded as an explicit version of a more conceptual one.
Consider a line bundle $L$ on $\P_k^r$ corresponding to $\mathcal{O}_{\P_k^r}(m)$ (such a line bundle is unique up to isomorphism) and choose a system of coordinates on the fibre $L_P$ for each $P \in \P^r(k)$. Then, $\ev_P (f)$ can be defined as the element of $k$ corresponding to $f_P \in L_P$ for this system of coordinates.
This is actually the genuine definition used by Manin to define algebraic--geometric codes in \cite{manin}.
Notice that another choice of coordinates on the fibres $L_P$ gives a Hamming--isometric code.
% In terms of invertible sheaves, choose for each point $P \in \P^r(k)$ a germ of section $s_P$ generating the stalk ${(\mathcal{O}_{\P^r}(m))}_{P}$ as a $\mathcal{O}_{\P^r\!,P}$-module. Then, for all $f\in \Gamma (\P^r , \mathcal{O}_{\P^r}(m))$ and all $P\in \P^r(k)$, there exists a unique function $\varphi_P$ which is regular at $P$ and such that $f=\varphi_P s_P$ and $\ev_P (f):=\varphi_{P}(P)$.
% In Definition \ref{evmap}, let $(X_0, \ldots , X_r)$ be a system of homogeneous coordinates on $\P^r$. Then, on the whole open set $U_0:=\{X_0 \neq 0\}$, the chosen generator of the stalk of a point is $X_0^m$.
% Thus, for $f\in \mathcal{F}_{m,r}$  the function $\varphi:=f/X_0^m$ is regular on $U_0$, satisfies $f= \varphi X_0^m$ and gives he evaluation map of Definition \ref{evmap} by setting $\ev_P (f):= \varphi (P)$.
\end{rem}

Now, let us define the notion of $m$--generality.

\begin{defn}[$m$--general position]\label{mgen}
 Let $m$ be a nonnegative integer.
 A family $P_1, \ldots ,$ $P_s$ of rational points of $\P^r$ is said to be in \textit{$m$--general position} if the evaluation maps $\ev_{P_1}, \ldots , \ev_{P_s}$ are linearly independent in $\Ff_{m,r}^{\vee}$.
\end{defn}

The following lemma gives a geometric interpretation for the notion of $m$--generality for $m\geq 1$.

\begin{lem}\label{lemmgen}
Let $m\geq 1$ be an integer and $P_1, \ldots , P_s$ be a set of rational points of $\P^r$. Then, the following assertions are equivalent.
 \begin{enumerate}[($i$)]
 \item\label{lgen1} The points $P_1, \ldots, P_s$ are in $m$--general position.
 \item\label{lgen2} For all $i \in \{1, \ldots , s\}$, there exists a hypersurface $H_i$ of degree $m$ in $\P^r$ containing the $P_j$'s for all $j\neq i$ and avoiding $P_i$.
 \item\label{lgen3} For all $i \in \{1, \ldots , s\}$, the point $P_i$ is not a base point of the linear system of hypersurfaces of degree $m$ in $\P^r$ containing all the $P_j$'s for $j\neq i$.
 \item\label{lgen4} The linear system $\Gamma$ of hypersurfaces of degree $m$ in $\P^r$ containing the points $P_1, \ldots , P_s$ has dimension
$$
\dim \Gamma= \dim \Ff_{m,r}-1-s.
$$
 \item\label{lgen5} $h^1 (\P^r, \mathcal{I}(m))=0$, where $\mathcal{I}$ is the ideal sheaf associated to the reduced zero--dimensional scheme supported by $P_1, \ldots, P_n$ and $\mathcal{I}(m)=\mathcal{I}\otimes \mathcal{O}(m)$.
 \end{enumerate}
\end{lem}

\begin{proof}
Proving (\ref{lgen1}) to (\ref{lgen4}) is an elementary exercise of linear algebra. For (\ref{lgen5}), consider the long exact sequence given by
$
0 \rightarrow \mathcal{I}(m) \rightarrow \mathcal{O}(m) \rightarrow \mathcal{S} \rightarrow 0,
$
where $\mathcal{S}$ is a skyscraper sheaf supported by $P_1, \ldots , P_n$.
\end{proof}

\begin{rem}
Notice that Definition \ref{mgen} makes sense even if $m=0$.
However, this case is removed in Lemma \ref{lemmgen} since items (\ref{lgen2}), (\ref{lgen3}) and (\ref{lgen4}) do not make sense for $m=0$.
\end{rem}

\begin{rem}\label{1gen}
  The notion of $1$--generality corresponds to the ``usual'' definition of general position, which is described at the beginning of the present section.
In $\P^r$, an $s$--tuple of points is in $1$--general position if the points are projectively independent, or equivalently if and only if they generate an $(s-1)$--dimensional linear subspace of $\P^r$.  
\end{rem}

\begin{defn}\label{mlink}
  A family $P_1, \ldots , P_s$ of rational points of $\P^r$ is said to be \textit{$m$--linked} if they are not in $m$--general position or equivalently if they fail to impose independent conditions on forms of degree $m$. It is said to be \textit{minimally $m$--linked} if it is $m$--linked and if each proper subset of $\{P_1, \ldots , P_s\}$ is in $m$--general position.
\end{defn}

We will see further that the notion of being \textit{minimally $m$--linked} is very useful for error--correcting codes.
Lemma \ref{remmin} gives some elementary algebraic and geometric translations of this definition which will be very often used in what follows. 

\begin{lem}\label{lemmin}
Let $m\geq 1$ be an integer.  A family $P_1, \ldots, P_s$ of rational points of $\P^r$ is minimally $m$--linked if and only if there exists a non-trivial relation of the form $\lambda_1 \ev_{P_1} + \cdots + \lambda_s \ev_{P_s}=0$ and that, for all such relation, the $\lambda_i$'s are all nonzero.
\end{lem}

\begin{proof}
It is an elementary exercise of linear algebra.
\end{proof}

\begin{rem}
  For dimensional reasons, one can prove easily that the number of elements of an $m$--general family of points in $\P^r$ is at most $\dim \Ff_{m,r}$ and that of a minimally $m$--linked family is at most $\dim \Ff_{m,r}+1$.
\end{rem}

\begin{rem}\label{remmin0}
Let $P_1, \ldots, P_s$ be a family of points in $\P^r$ and let $m$ be a nonnegative integer. Assume that $s\leq \dim \Ff_{m,r}$, then the $P_i$'s are minimally $m$--linked if and only if for all $i_0\in \{1, \ldots, s\}$ the linear system of hypersurfaces of degree $m$ containing the points $P_1, \ldots, P_{i_0-1},P_{i_0+1}, \ldots , P_s$ is nonempty and has $P_{i_0}$ as a base point.
\end{rem}

\begin{rem}\label{remmin}
The previous remark entails that, to prove that a family of points $P_1, \ldots,$ $P_s\in \P^r$ with $s\leq \dim \Ff_{m,r}$ is \textit{not} minimally $m$--linked, it is sufficient to prove that for one of these points $P_{i_0}$, there exists a hypersurface of degree $m$ containing the $P_j$'s for $j\neq i_0$ and avoiding $P_{i_0}$.
\end{rem}

We conclude the present section with Lemma \ref{geom}, which is crucial in the present paper.
Indeed, it enables to work over an algebraically closed field of the form $\overline{\F}_q$ in order to get information on the minimum distance of some codes, even if a code is a vector space of a finite field $\F_q$.
Such a ``geometrisation'' of the problem is very useful since over infinite fields, the positive--dimensional linear systems have infinitely many elements.

\begin{lem}\label{geom}
Let $P_1, \ldots , P_s$ be a family of $k$--rational points of $\P^r$. Let $L$ be an algebraic extension of $k$. Then, the points $P_1, \ldots, P_s$ are in $m$--general position (resp. are $m$--linked, resp. are minimally $m$--linked) in $\P^r_k$ if and only if they are in $m$--general position (resp. are $m$--linked, resp. are minimally $m$--linked) in $\P^r_L$.
\end{lem}

\begin{proof}
  Linearly independent (resp. linked) vectors in $\Ff_{m,r}(k)^{\vee}$ remain independent (resp. linked) as vectors of $\Ff_{m,r}(L)^{\vee}=\Ff_{m,r}(k)^{\vee}\otimes_k L$.
\end{proof}

\section{Duals of algebraic--geometric codes}\label{secdual}

In what follows, when we deal with algebraic--geometric codes and only in this situation (that is in the present section and in Section \ref{secex}), we always stay in the following context.

\subsection{Context and notations}\label{codescontext}
In what follows, $X$ is a smooth geometrically connected projective variety over $\F_q$, which is a \textbf{complete intersection} in some projective space $\P^r$ for some $r\geq 2$.
Moreover, $m$ is a nonnegative integer and $G_m$ is a divisor on $X$ which is linearly equivalent to a scheme--theoretic intersection of $X$ with a hypersurface of degree $m$.
In addition, $P_1, \ldots, P_n$ is a family of rational points of $X$ avoiding the support of $G_m$ and we denote by $\Delta$ the $0$--cycle $\Delta:=P_1+ \cdots +P_n$.

 From \cite[Exercise II.8.4]{H}, the variety $X$ is \textit{projectively normal} (see \cite[Exercise I.3.18]{H} for a definition) and an element of $L(G_m)$ can be identified to a restriction to $X$ of an element of $\Ff_{m,r}$.
The connection between minimum distance of $C_L (\Delta, G_m)^{\bot}$ and the notion of $m$--generality lies in the elementary Lemma \ref{codeword} below.

\subsection{Codewords of the dual and configurations of points} 
First, let us notice a usual abuse of language, in the next sections.

\medbreak

\noindent\textbf{Abuse of language.} In what follows, given a codeword $c \in C_L (\Delta, G_m)$ or $c \in C_L (\Delta, G_m)^{\bot}$, we will call \textit{support of $c$} the set of points $P_{i_1}, \ldots, P_{i_s}$ in $\supp (\Delta)$ corresponding to the nonzero coordinates of $c$.

Thanks to the following proposition, the problem of finding a lower bound for the minimum distance of the code $C_L (\Delta, G_m)^{\bot}$ is translated into that of  finding configurations of (minimally) $m$--linked points in the support of $\Delta$.

\begin{prop}\label{dminlink}
  The minimum distance of the code $C_L (\Delta, G_m)^{\bot}$ is the smallest number of $m$--linked points in the support of $\Delta$.
\end{prop}

\begin{rem}
  Equivalently it is the smallest number of \textit{minimally} $m$--linked points of $\supp (\Delta)$.
\end{rem}

The proof of Proposition \ref{dminlink} is a straightforward consequence of the  following lemma.

\begin{lem}\label{codeword}
  There exists a nonzero codeword $c\in C_L (\Delta, G_m)^{\bot}$ with support contained in $\{P_{i_1}, \ldots , P_{i_s}\}$ if and only if these points are $m$--linked.
Furthermore, if these points are minimally $m$--linked, then the support of such a codeword is equal to $\{P_{i_1}, \ldots,P_{i_s}\}$.
\end{lem}

\begin{proof}
The existence of the codeword $c\in C_L (\Delta, G_m)^{\bot}$ with support $\{P_{i_1}, \ldots, P_{i_s}\}$ entails that of a nonzero linear relation linking the evaluation maps $\ev_{P_{i_1}}, \ldots , \ev_{P_{i_s}}$ in $\Ff_{m,r}^{\vee}$. Conversely, if $P_{i_1}, \ldots, P_{i_s}$ are $m$-linked, then a non-trivial linear relation linking the corresponding evaluation maps entails the existence of a nonzero codeword with support contained in $\{P_{i_1}, \ldots, P_{i_s}\}$.
If the points are minimally $m$--linked, then, from Lemma \ref{lemmin}, a non-trivial linear relation linking the corresponding evaluation maps gives a codeword with support equal to $\{P_{i_1}, \ldots, P_{i_s}\}$.
\end{proof}

Therefore, minimally $m$--linked configurations of points seem to be useful to estimate the minimum distance of $C_L (\Delta, G_m)^{\bot}$.
Let us state the main results concerning the minimum distance of the dual of a functional code.

\subsection{Lower bounds for the minimum distance of the dual code}
To state some results on the minimum distance of the codes of the form $C_L (\Delta, G_m)^{\bot}$, we will treat separately the ``small'' values of $m$, i.e. $m=0$ and $1$ and the other ones, i.e. $m\geq 2$.

\subsubsection{Small values of $m$}\label{small}

If $m=0$, then 
the code $C_L(\Delta,G_0)$ is a \textit{pure repetition code}, i.e. it is generated over $\F_q$ by the codeword $(1, \ldots ,1)$.
Thus, the minimum distance of $C_L (\Delta,G_0)^{\bot}$ is $2$ and any pair of distinct points $P_i,P_j\in \supp (\Delta)$ is the support of a codeword in $C_L (\Delta,G_0)^{\bot}$.
In terms of $m$-generality one sees easily that one point is always $0$--general and that two distinct points are always $0$--linked and hence minimally $0$--linked.

\medbreak

If $m=1$, then we have the following result.

\begin{lem}\label{meq1}
  In the context described in Section\ref{codescontext}, if for all $t\leq n-2$, no $t+2$ of the $P_i$'s lie on a linear subspace of dimension $t$, then the minimum distance of $C_L(\Delta, G_1)^{\bot}$ is $n$.
Moreover, let $s$ be the smallest integer such that there exist $s+2$ elements of $\supp (\Delta)$ lying in a linear subspace of dimension $s$. Then $s+2$ is the minimum distance of the code $C_L (\Delta, G_1)^{\bot}$. 
\end{lem}

\begin{proof}
From Remark \ref{1gen}, a $t$--tuple of points of $\P^r$ is $1$--general if and only if it generates a linear subspace of dimension $t-1$.
If the integer $s$ of the statement exists, then the smallest number of $1$--linked points of $\supp({\Delta})$ is $s+2$ and, from Proposition \ref{dminlink}, this gives the minimum distance of $C_L(\Delta,G_1)^{\bot}$.
If $s$ does not exist, then the minimum distance of $C_L(\Delta,G_1)^{\bot}$ is obviously $n$.
\end{proof}

\subsubsection{Other values of $m$}

\begin{thm}\label{main}
  In the context described in Section \ref{codescontext}, let $m$ be an integer greater than or equal to $2$ and $d$ be the minimum distance of the code $C_L (\Delta, G_m)^{\bot}$. Then,
  \begin{enumerate}[(i)]
  \item\label{T1} $d=m+2$ if and only if $m+2$ of the $P_i$'s are collinear in $\P^r$;
 \item\label{T2} $d=2m+2$ if and only if no $m+2$ of the $P_i$'s are collinear and $2m+2$ of them lie on a plane conic (possibly reducible);
 \item\label{T3} $d=3m$ if and only if no $m+2$ of the $P_i$'s are collinear, no $2m+2$ of them lie on a plane conic and $3m$ of them are coplanar and lie at the intersection of a cubic and a curve of degree $m$ having no common irreducible component;
 \item\label{T4} $d\geq 3m+1$ if and only if no sub-family of the $P_i$'s satisfies one of the three above-cited configurations.
  \end{enumerate}
Moreover, in case (\ref{T1}) (resp. (\ref{T2}), resp. (\ref{T3})), the minimum weight codewords are supported by the configuration of points in question.
\end{thm}

\begin{rem}\label{meq2}
  If $m=2$, then the condition of Theorem \ref{main}(\ref{T3}) cannot happen.
%   Indeed, on the one hand, one expects that no $6$ of the $P_i$'s lie on a conic, on the other hand, one expects that $6$ of th $P_i$'s lie at the intersection of a cubic and a curve of degree $2$, hence on a conic.
Consequently, in this situation, the statement is simplified: \textit{the minimum distance $d$ of $C_L (\Delta, G_2)^{\bot}$ is
  \begin{enumerate}[(i)]
  \item $4$ if and only if $4$ of the $P_i$'s are collinear;
  \item $6$ if and only if $6$ of the $P_i$'s lie on a plane conic;
  \item $\geq 7$ if and only if none of the above-cited configurations happens.
  \end{enumerate}
}
\noindent Therefore, in Section \ref{secfourth}, which is devoted to the end of the proof of Theorem \ref{main}, we assume that $m\geq 3$.
\end{rem}

\begin{rem}
  If $m\geq 2$, one checks that $m+2$, $2m+2$ and $3m$ are lower than or equal to $\dim \Ff_{m,r}$ (recall that we are in the context of Section \ref{codescontext} and hence $r\geq 2$). Therefore, to prove that $m+2$, $2m+2$ or $3m$ points are (resp. are not) minimally $m$--linked, one can use Remark \ref{remmin0} (resp. Remark \ref{remmin}).  
\end{rem}

%----------------------------------------------------------------------------

% \noindent \hrulefill Peut \^etre coup\'e \hrulefill

% Theorem \ref{main} gives estimates for the minimum distance of the dual of an algebraic geometric code.
% Moreover if the lower bound is reached, then it gives a geometric description for the minimum weight codewords.

% Actually, the notions of being $m$--linked can help to describe the supports of all the codewords of weight $\leq 3m$.
% This is the purpose of the next statement.

% \begin{thm}[Small weight codewords of the code $C_L (\Delta, G_m)^{\bot}$]\label{weights}
% Under the conditions of Theorem \ref{main}, let $c$ be a codeword of $C_L(\Delta, G_m)^{\bot}$. Then
% \begin{enumerate}[(i)]
% \item\label{W1} if $m+2 \leq w(c) \leq 2m+1$, then the support of $c$ consists in collinear points;
% \item\label{W2} if $2m+2 \leq w(c) \leq 3m-1$, then the support of $c$ consists in coplanar points lying on a conic,
% \end{enumerate}
% where $w(c)$ denotes the Hamming weight of $c$.
% \end{thm}

% \noindent \hrulefill Fin de coupe \hrulefill

%---------------------------------------------------------------------------

To prove Theorem \ref{main}, we will actually prove the following statement, which is a ``geometric version'' of Theorem \ref{main}. 

\begin{thm}[Geometric version of Theorem \ref{main}]\label{geomain}
Let $P_1, \ldots , P_n$ be a family of distinct points in a projective space $\P^N$ and let $m \geq 2$ be an integer.
Then, the smallest number of $m$--linked points in $\{P_1, \ldots , P_n\}$ is
\begin{enumerate}[(i)]
\item\label{G1} $m+2$ if and only if $m+2$ of the $P_i$'s are collinear;
\item\label{G2} $2m+2$ if and only if no $m+2$ of the $P_i$'s are collinear and $2m+2$ of the $P_i$'s lie on a plane conic;
\item\label{G3} $3m$ if and only if no $m+2$ of the $P_i$'s are collinear, no $2m+2$ of them lie on a plane conic and $3m$ of them lie at the intersection of two coplanar plane curves of respective degrees $3$ and $m$;
\item\label{G4} $>3m$ if and only if the $P_i$'s do not satisfy any of the above configurations.
\end{enumerate}
\end{thm}

The proof of Theorem \ref{geomain} will be the purpose of Sections \ref{secfirst} to \ref{secfourth}.
The organisation of this proof is detailed in Section \ref{detail} below.
First let us show that Theorem \ref{geomain} entails Theorem \ref{main}.

\begin{proof}[Proof of Theorem \ref{geomain} $\Rightarrow$ Theorem \ref{main}]
  Proposition \ref{dminlink} asserts that the minimum distance of $C_L (\Delta, G_m)^{\bot}$ equals the smallest number of the $P_i$'s which are $m$--linked (and hence minimally $m$--linked). Therefore, Theorem \ref{geomain}(X) $\Rightarrow$ Theorem \ref{main}(X) for all $X$ in $\{$\ref{T1}, \ref{T2}, \ref{T3}, \ref{T4}$\}$.
\end{proof}

\subsection{The proof of Theorem \ref{geomain}}\label{detail}

It is worth noting that Cayley--Bacharach theorem (see Theorem \ref{intersec} further) asserts that the configurations described in Theorem \ref{geomain} are minimally $m$--linked and hence, from Lemma \ref{codeword}, provide supports of codewords in $C_L (\Delta, G_m)^{\bot}$. The point of the proof is to make sure that these minimally $m$--linked configurations are the smallest ones. In particular, an interesting step of the proof is Lemma \ref{coplan1} which asserts that, whatever the ambient dimension is, a minimally $m$--linked configuration of less than $3m$ points is contained in a $2$--dimensional linear subspace.

 In Section \ref{secmlink}, Lemma \ref{indmin} and Theorem \ref{intersec} are stated.
Lemma \ref{indmin} is a nice trick to handle the notion of being \textit{minimally $m$--linked}.
Theorem \ref{intersec} is one of the numerous versions of Cayley--Bacharach and gives the description plenty of configurations of minimally $m$--linked points.
Among others things, it asserts that $m+2$ collinear points, $2m+2$ points on a plane conic or $3m$ points lying at the intersection of two coplanar curves with respective degrees $m$ and $3$ are minimally $m$--linked.

\medbreak

In Section \ref{secfirst}, one proves Proposition \ref{mg1} asserting that less than $m+1$ points of $\P^r$ are always in $m$--general position.
Proposition \ref{mg1} together with Theorem \ref{intersec} (applied to $a=1$) entails obviously the \textit{``if''} part of Theorem \ref{geomain}(\ref{G1}).
Conversely, we prove Proposition \ref{convmg1} which asserts that, any $m+2$ points which are $m$--linked are collinear.
This yields the \textit{``only if''} part of Theorem \ref{geomain}(\ref{G1}).

\medbreak

 In Section \ref{secsecond}, we prove Proposition \ref{mg2}, which asserts that any set of at most $2m+1$ points of $\P^r$ such that no $m+2$ of them are collinear is in $m$--general position.
Proposition \ref{mg2} in addition with Theorem \ref{intersec} (applied to $a=2$), entails the \textit{``if''} part of Theorem \ref{geomain}(\ref{G2}).
Conversely, one proves Proposition \ref{convmg2}, which asserts that, any $m$--linked configuration of $2m+2$ points such that no $m+2$ of them are collinear lies on a plane conic.
This yields the \textit{``only if''} part of Theorem \ref{geomain}(\ref{G2}).

%-----------------------------------------------------------------------------
% \noindent \hrulefill Peut être coupé \hrulefill

% One concludes Section \ref{secsecond} with the proof of Theorem \ref{weights}(\ref{W1}).

% \noindent \hrulefill Fin de coupe \hrulefill
%-----------------------------------------------------------------------------

\medbreak

 In Section \ref{secthird}, we prove Proposition \ref{mg3}, which asserts that any set of at most $3m-1$ points of $\P^r$ such that no $m+2$ of them are collinear and no $2m+2$ of them lie on a plane conic, is in $m$--general position.
Proposition \ref{mg3} in addition with Theorem \ref{intersec} (applied to $a=3$) yields the \textit{``if''} part of Theorem \ref{geomain}(\ref{G3}).

%-----------------------------------------------------------------------------
% \noindent \hrulefill Peut être coupé \hrulefill

% One concludes Section \ref{secthird} with the proof of Theorem \ref{weights}(\ref{W2}).

% \noindent \hrulefill Fin de coupe \hrulefill
%-----------------------------------------------------------------------------

\medbreak 

Section \ref{secfourth} is devoted to the proof of Proposition \ref{convmg3}, which asserts that any $m$--linked configuration of $3m$ points such that no $m+2$ of them are collinear and no $2m+2$ of them are on a plane conic is a set of coplanar points lying at the intersection of a cubic and a curve of degree $m$ having no common component.
This concludes the proof of Theorem \ref{geomain} since it yields the \textit{``only if''} part of (\ref{G3}) and (\ref{G4}).

\medbreak

Before starting the different steps of the proof of Theorem \ref{main}, let us present some applications of it.

\section{Examples and applications}\label{secex}

Even if the objective of the present article is to get results on duals of algebraic--geometric codes on higher--dimensional varieties, Theorem \ref{main} holds for varieties of any dimension.
Surprisingly, when the variety $X$ is a plane curve, Theorem \ref{main} gives a relevant lower bound for the minimum distance of some algebraic--geometric codes on $X$.

\subsection{Algebraic--geometric codes on plane curves}\label{curves}

\subsubsection{Context}\label{curvecont}
Let $a$ be a positive integer.
Let $X\subset \P^2$ be a smooth projective plane curve of degree $a$ over $\F_q$.
Let $m$ be a nonnegative integer, $L$ be a line of $\P^2$ and $G_m$ be the pullback of $mL$ by the natural inclusion map $X \hookrightarrow \P^2$.
Let $P_1, \ldots , P_n$ be $n$ rational points of $X$ avoiding the support of $G_m$ and denote by $D$ the divisor $D:=P_1+ \cdots +P_n$.

\subsubsection{The code $C_L (D,G_m)^{\bot}$} From \cite[Theorem II.2.8]{sti}, the dual $C_L (D,G_m)^{\bot}$ of the functional code is the \textit{differential} code denoted by $C_{\Omega}(D, G_m)$.
Denote by $d$ the minimum distance of $C_L (D,G_m)^{\bot}$.
Let $\delta_G$ be the \textit{Goppa designed distance}. From \cite[Theorem II.2.7]{sti}, we have
$
\delta_G=\deg(G_m)-(2g_X-2),
$
where $g_X$ denotes the genus of $X$, which is
$
g_X =(a-1)(a-2)/2
$.
This gives
\begin{equation}
  \label{deltaG}
  \delta_G=a(m+3-a).
\end{equation}

We know that $d\geq \delta_G$.
Let us study the lower bound for $d$ given by Theorem \ref{main}.

\subsubsection{Lower bound for the dual minimum distance}\label{lowerb}
First, notice that if the degree of the curve $X$ is $1$ or $2$, then $X$ is isomorphic to $\P^1$ and the codes on it are Reed--Solomon codes, for which the Goppa designed distance equals the genuine distance (which reaches the Singleton bound) and hence is optimal.
Therefore, from now on, assume that the degree $a$ of $X$ is greater than or equal to $3$.

Denote by $\delta$ the lower bound for the minimum distance given by Theorem \ref{main}:
$$
\delta=\left\{
  \begin{array}{ccccccc}
    m+2 & \textrm{if} & 0 & \leq & m & \leq & a-2 \\
    2m+2 & \textrm{if} &  &  & m & = & a-1 \\
    3m & \textrm{if} &  &  & m & \geq & a \\
  \end{array}
\right. .
$$

\noindent Notice that $\delta$ is always positive, which is not true for the Goppa designed distance $\delta_G$.
Therefore, $\delta$ gives a relevant lower bound for the minimum distance of $C_{L}(D, G_m)^{\bot}$ when $\delta_G \leq 0$.

\begin{thm}[Minimum distance for codes on curves]\label{Goppa}
  Let $X$ be a smooth plane curve of degree $a\geq 3$ and consider the code $C_{\Omega}(D,G_m)=C_L(D,G_m)^{\bot}$.
Then, $\delta > \delta_G$ if and only if $\delta_G \leq 0$ (or equivalently, if and only if $m \leq a-3$). 
In other words, $\delta$ improves the Goppa designed distance $\delta_G$ as a lower bound for the minimum distance of the code whenever $\delta_G$ is negative and hence irrelevant for coding theory.
\end{thm}

\begin{proof}
Let us compare the numbers $\delta$ and $\delta_G$.
Using (\ref{deltaG}), a brief computation gives
$$ 
\delta- \delta_G=
\left\{
  \begin{array}{ccccccc}
    (a-1)(a-2-m) & \textrm{if} & 0 & \leq & m & \leq & a-2 \\
    0 & \textrm{if} &  &  & m & = & a-1 \\
    (a-3)(a-m) & \textrm{if} &  &  & m & \geq & a \\
  \end{array}
\right. .
$$

\noindent Consequently, $\delta - \delta_G > 0$ if and only if $m \leq a-2$.
That is, from (\ref{deltaG}), this difference is nonnegative if and only if the Goppa designed distance $\delta_G$ is negative.
\end{proof}

\begin{rem}
In the proof, one can also see that $\delta = \delta_G$ for all $m \in \{a-3,a\}$.
\end{rem}

\begin{exmp}\label{ex}
  Consider the finite field $\F_{64}$ and the curve $C$ of equation
$$
F_C:=w^{24} x^{11} + w^{44} x^6 y^2 z^3 + w^{24} x^5 y z^5 + w^{20} x^4 y^6 z + w^{33} x^2 z^9 +
$$
$$
     w^{46} x y^5 z^5 + w^{46} x z^{10} + w^{39} y^{11} + w^{30} y^2 z^9,
$$
where $w$ is a primitive element of $\F_{64}$ over $\F_2$ with minimal polynomial $x^6 + x^4 + x^3 + x + 1$.
This curve has $80$ rational points in the affine chart $\{z \neq 0\}$ and $1$ rational point \textit{at infinity}.
Using the previous results, one sees that the Goppa designed distance of $C_L(D, G_m)^{\bot}$ is negative for $m\leq 8$.
Using Theorem \ref{main}, we prove that the codes $C_L(D,G_m)^{\bot}$ for $m=1, \ldots, 8$ are respectively of the form:
$[80,77,\geq 3]$, $[80,74,\geq 4],$ $[80,70,\geq 5],$ $[80,65,\geq 6]$, $[80,59,\geq 7]$,
$[80,52,\geq 8]$, $[80,46,\geq 9]$ and $[80,35,\geq 10]$.
\end{exmp}

Afterwards, under some geometric condition on the points $P_1, \ldots, P_n$, one can improve the Goppa designed distance by using Theorem \ref{main}.
It is worth noting that if the lower bound $m+2$ is not reached (that is, if no $m+2$ of the $P_i$'s are collinear), then this bound jumps directly to $2m+2$.
By this way one can get, under some non--incidence conditions, some good improvements of the Goppa bound even if it is positive.

\begin{thm}\label{Improve}
  Under the assumptions of Theorem \ref{Goppa}, 
  \begin{enumerate}
  \item if $m \leq a-2$ and no $m+2$ of the $P_i$'s are collinear, then the minimum distance $d$ of $C_L (D, G_m)^{\bot}$ satisfies $d\geq 2m+2$ and this bound improves that of Goppa;
  \item if $m \leq a-1$, no $m+2$ of the $P_i$'s are collinear and no $2m+2$ of them lie on a conic,
    then $d \geq 3m$ and this bound improves that of Goppa;
  \item if $m \leq a$, the $P_i$'s do not satisfy any of the above condition and no $3m$ of them lie on a cubic, then
$d \geq 3m+1$ and this bound improves that of Goppa.
  \end{enumerate}
\end{thm}

  \begin{proof}
    It is a straightforward consequence of Theorem \ref{main}.
  \end{proof}

\medbreak

\begin{exmp}
  Back to Example \ref{ex}, a computation using the software \textsc{Magma} yields only one line containing at least 7 of the $P_i$'s. It is the line $L$ of equation $x=0$, which contains $10$ of the $P_i$'s.
Therefore by removing $4$ (resp. $3$, resp. $2$, resp. $1$) of the $P_i$'s on $L$, one gets a divisor $D^{(4)}$ (resp. $D^{(3)}$, resp. $D^{(2)}$, resp. $D^{(1)}$) and the codes $C_L (D^{(i)}, G_{4+i})^{\bot}$ for $i\in \{1,2,3,4\}$ are respectively of the form
$[76, 55, \geq 12]$, $[77,49, \geq 14]$, $[78,44, \geq 16]$ and $[79,34, \geq 18]$.

Moreover, the Goppa designed distance asserts that $C_L (D, G_{9})^{\bot}$ has a minimum distance greater than or equal to $11$.
However, since no $11$ of the $P_i$'s are collinear, Theorem \ref{Improve} asserts that this minimum distance is greater than or equal to $20$.
Thus, \textbf{the obtained lower bound is $9$ units bigger than that of Goppa}.
\end{exmp}

The previous example presents actually a good method to get good codes on curves by selecting the points of evaluation
.
Indeed, assume there are only few lines (resp. conics, resp. cubics) containing $m+2$ (resp. $2m+2$, resp. $3m$) of the $P_i$'s.
Then one can remove some points of these lines (resp. conics, resp. cubics) such that the lower bound for the minimum distance jumps to $2m+2$ (resp. $3m$, resp. $3m+1$).

Further, in Section \ref{back}, we give an interpretation of the Goppa designed distance for plane curves in terms of minimally $m$-linked points in $\P^2$.
% and
% \begin{enumerate}
% \item a geometric necessary and sufficient condition for the Goppa designed distance to be reached;
% \item a geometric description of the minimum weight codewords.
% \end{enumerate}

\subsection{Surfaces in $\P^3$}\label{surfaces}

Here, we assume that $q\geq 3$.
The binary case will be treated in Section \ref{binary}.

\subsubsection{Context} Let $a$ be a positive integer and $X$ be a smooth projective geometrically connected surface of degree $a$ defined over $\F_q$ and embedded in $\P^3$.
Let $H$ be a plane of $\P^3$, let $m$ be a nonnegative integer and $G_m$ be the pullback of the divisor $mH$ by the canonical inclusion $X \hookrightarrow \P^3$.
Let $P_1,  \ldots , P_n$ be a family of rational points of $X$ avoiding the support of $G_m$ and $\Delta$ be the $0$--cycle $\Delta:=P_1+ \cdots +P_n$.

% \subsubsection{The binary case}
% Let us first consider the case when the base field is $\F_2$, for which the lower bounds for the minimum distance of the dual code does not depend on the choice of the surface.

% \begin{thm}[The binary Case]
%   When the base field is $\F_2$, the minimum distance of $C_L (\Delta, G_m)^{\bot}$ satisfies
% $
% \left\{
%   \begin{array}{lcl}
%     d=4 & \textrm{if} & m =1\\ 
%     d> 3m & \textrm{if} & m> 1
%   \end{array}
% \right.
% $, if $X$ is either elliptic or hyperbolic.
% \end{thm}

% \begin{proof}
%   Just notice that, since the $P_i$'s are contained in an affine chart of $\P^3$, a plane section of $X$ cannot contain more than $\sharp \mathbf{A}^2(\F_2)=4$ of the $P_i$'s. The result is a straightforward consequence of this assertion an Theorem \ref{main}.
% \end{proof}

\subsubsection{Duals of codes on quadrics}
Let $X$ be a quadric in $\P^3$.
There are two isomorphism classes of smooth quadrics in $\P^3$, respectively called \textit{hyperbolic} and \textit{elliptic} quadrics.
Hyperbolic quadrics contain lines defined over $\F_q$ and elliptic quadrics do not.

For each isomorphism class, there exists an affine chart $U$ of $X$ containing exactly $q^2$ rational points. One chooses the complement of $U$ to be the support of $G_m$ and the sum of the rational points of $U$ to be $\Delta$.

\begin{thm}
The minimum distance $d$ of $C_L(\Delta, G_m)^{\bot}$ satisfies the following relations.

If $X$is hyperbolic, then $\left\{
  \begin{array}{lcl}
    d=m+2 & \textrm{if} & m\leq q-2\\
    d=2m+2 & \textrm{if} & m=q-1\\
    d= 3m & \textrm{if} & m= q \\
    d>3m & \textrm{if} & m>q
  \end{array}
\right.$.

\smallbreak

If $X$ is elliptic, then $
\left\{
  \begin{array}{lcl}
    d=2m+2 & \textrm{if} & m \leq (q-1)/2\\ 
    d> 3m & \textrm{if} & m> (q-1)/2
  \end{array}
\right.
$.
\end{thm}

\begin{proof}
Notice that, since the $P_i$'s all lie on an affine chart of $\P^3$, no $q+1$ of them are collinear and no $2q+1$ of them lie on a conic.
If $X$ is hyperbolic (resp. elliptic), then plane sections of $X$ are either irreducible plane conics containing at most $q+1$ of the $P_i$'s, or a union of two rational lines (resp. a union of two lines defined over $\F_{q^2}$ and conjugated by the Frobenius) containing at most $2q$ of the $P_i$'s (resp. containing $1$ of the $P_i$'s).
This description of the plane sections of $X$ together with Theorem \ref{main} leads easily to the expected result.
\end{proof}

\begin{exmp}
For $q=3$ one gets codes of the following form.
$$
\begin{array}{|c|c|c|}
  \hline
   & X\ \textrm{is hyperbolic} & X\ \textrm{is elliptic} \\
  \hline
   m=1 &[9,5,3] & [9,5,4] \\
  \hline
\end{array}
$$
For $q=4$ one gets codes of the following form.
$$
\begin{array}{|c|c|c|}
  \hline
   & X\ \textrm{is hyperbolic} & X\ \textrm{is elliptic} \\
  \hline
   m=1 &[16,12,3] & [16,12,4] \\
  \hline
   m=2 & [16,7,4] & [16,7,\geq 7]\\
  \hline
\end{array}
$$
For $q=5$ one gets codes of the following form.
$$
\begin{array}{|c|c|c|}
  \hline
   & X\ \textrm{is hyperbolic} & X\ \textrm{is elliptic} \\
  \hline
   m=1 &[25,21,3] & [25,25,4] \\
  \hline
   m=2 & [25,16,4] & [25,16,6]\\
  \hline
   m=3 & [25,9,5] & [25,9, \geq 9]\\
  \hline
\end{array}
$$
\end{exmp}

\subsubsection{Duals of codes on cubics}
Let $X$ be a cubic in $\P^3$.
As in the previous case we state a result by separating the cases when $X$ contains rational lines and when it does not.
Indeed, even if a cubic surface always contains $27$ lines over the algebraic closure of the base field, all these lines can be non--rational (see \cite[Chapter 3]{kollar}).

\begin{thm}\label{cubique}
  The minimum distance $d$ of $C_L (\Delta , G_m)^{\bot}$ satisfies the following relations.

If $X$ contains rational lines, then $
 \left\{
  \begin{array}{lcl}
    d=m+2 & \textrm{if} & m\leq q-2\\
    d=2m+2 & \textrm{if} & m=q-1\\
    d\geq 3m & \textrm{if} & m\geq q 
  \end{array}
\right.
$.

If $X$ does not contain any rational line, then $
 \left\{
  \begin{array}{lcl}
    d\geq 3 & \textrm{if} & m= 1\\
    d\geq 6 & \textrm{if} & m= 2\\
    d\geq 3m & \textrm{if} & m\geq q 
  \end{array}
\right.
$.
\end{thm}

\begin{proof}
As in the previous example, since the $P_i$'s lie on an affine chart of $\P^3$, no $q+1$ of them are collinear and no $2q+1$ of them lie on a conic.
Moreover, if the cubic surface $X$ does not contain rational lines, then it does not contain any rational plane conic.
This yields the result thanks to Theorem \ref{main}.  
\end{proof}

%Some of these results are proved in \cite{oim2} Section 8 by a different manner.

\medbreak

\begin{exmp}
  The Hermitian surface over $\F_4$ is the surface of equation 
$
x^3+y^3+z^3+t^3=0
$.
This surface has $36$ rational points in the affine chart $\{t\neq 0\}$ and contains plenty of lines.
The code $C_L (\Delta, G_1)^{\bot}$ is $[36,32,3]$ and the supports of the codewords of weight $3$ are triples of collinear points.
The code $C_L (\Delta, G_2)^{\bot}$ is $[36,26,4]$ and the supports of the codewords of weight $3$ are $4$--tuples of collinear points.
The code $C_L (\Delta, G_3)^{\bot}$ is $[36,17,8]$ and the supports of the codewords of weight $3$ are $8$--tuples of points lying on plane conics (since $q=4$, such conics are reducible).
\end{exmp}

\begin{exmp}
  In \cite{zarzar}, an example of a cubic surface over $\F_9$ containing no rational lines is given.
The author proves that on this surface, the code $C_L (\Delta, G_2)$ is a $[100,10,68]$ code.
Using Theorem \ref{cubique}, one proves that its dual is a $[100,90, \geq 6]$ code.
Theorem \ref{cubique} asserts also that $C_L (\Delta, G_3)^{\bot}$ is $[100,81, \geq 9]$.
\end{exmp}

\begin{exmp}
  Another example is given in \cite{agctvoloch}: the surface over $\F_3$ defined by the affine equation $x^3+y^3+z^3-zx^2-yx^2-yz^2+xz^2+1$.
The code $C_L (\Delta, G_1)$ on this surface is $[13,4,7]$.
From Theorem \ref{cubique}, its dual is a $[13,9,\geq 3]$ code.

Moreover, the authors also assert that this surface does not contain any rational line over $\F_9$.
They prove that, over $\F_9$, the code $C_L (\Delta, G_2)$ is $[91,10,61]$. Theorem \ref{cubique} entails that its dual is $[91,81,\geq 6]$.
Moreover, Theorem \ref{cubique} entails that $C_L(\Delta, G_3)^{\bot}$ is a $[91,72, \geq 9]$ code over $\F_9$.
\end{exmp}

\subsubsection{Surfaces of higher degree}
To conclude this subsection on codes on surfaces, let us give some example of surfaces of higher degree.
Theorem \ref{inclusion} together with Remark \ref{degree4} in Appendix \ref{append} asserts that almost all surfaces in $\P^3$ of degree $\geq 4$ do not contain any line, plane conic and plane cubic, even over the algebraic closure of their based field.
Moreover, we produced a \textsc{Magma} program checking all the plane sections of a surface and asserting whether they are all irreducible.

Thus, one can expect to find a lot of surfaces giving dual codes of minimum distance $>3m$.

\begin{exmp}
  Over $\F_7$, the surface defined by the equation
$x^4 + 2 x^3 y + 4 x^3 t + 3 x^2 z^2 + 6 x y^3 + 4 x z^2 t + 4 y^3 z + 6 y^2 t^2 
    + 5 y t^3 + 4 z^4
 $
does not contain any line, plane conic or plane cubic.
It has $54$ rational points in the affine chart $\{t \neq 0\}$. Therefore, Theorem \ref{main} asserts that the codes $C_L(\Delta, G_{m})^{\bot}$ are respectively of the form $[54,50,\geq 3]$, $[54, 44, \geq 4]$, $[54, 34, \geq 9]$ and $[54,20,\geq 12]$ when $m=1,2,3,4$.
\end{exmp}

\begin{exmp}
  Over $\F_8$, the surface defined by the equation $\gamma^2 x^5 + x^4 y + \gamma^5 x^4 z + \gamma^4 x^3 z t + \gamma^6 x^2 z^3 + \gamma^4 x y t^3 + 
    \gamma^3 x z^4 + \gamma^5 y^4 t + \gamma^3 y^2 t^3 + \gamma^6 y z^4 + \gamma^5 y t^4 + \gamma^5 z^2 t^3$, where $\gamma$ denotes a primitive element of $\F_8/\F_2$, contains also no line, plane conic or plane cubic.
Its affine chart $\{t\neq 0\}$ contains $64$ rational points and hence, the codes $C_L(\Delta, G_{m})^{\bot}$ are of the form 
$[64,60,\geq 3]$, $[64, 54, \geq 4]$, $[64, 44, \geq 5], [64,29,\geq 12], [64,9, \geq 15]$, when $m=1,\ldots , 5$.

\end{exmp}

\subsection{Higher--dimensional varieties}\label{ldpc}

For higher--dimensional varieties, the situation is more difficult, since it is quite harder to check whether a variety contains a line (resp. a plane conic) or not.

However, Theorem \ref{inclusion} in Appendix \ref{append} gives some \textit{generic} results on codes on hypersurfaces of fixed degree.

For instance, it asserts that in $\P^4$, almost all hypersurfaces of degree $a\geq 6$ do not contain any line, plane conic or plane cubic.
Therefore, we know that codes $C_L (\Delta, G_m)^{\bot}$ have minimum distance $d \left\{
  \begin{array}{ccl}
   \geq m+2 & \textrm{if} & m \leq a-2\\
   \geq  2m+2 & \textrm{if} & m= a-1\\
   \geq 3m & \textrm{if} & m=a \\
   > 3m & \textrm{if} & m>a
  \end{array}
\right.$.

\subsection{Binary codes}\label{binary}
To conclude this section let us consider the case of algebraic-geometric codes over $\F_2$.

\begin{thm}
  Let $H$ be a hypersurface of $\P^N_{\F_2}$ with $N\geq 3$, let $G_m$ be $m$ times a hyperplane section of $H$ and $\Delta$ be a formal sum of points avoiding the support of $G_m$. Then the minimum distance $d$ of $C_L(\Delta, G_m)^{\bot}$
is $\left\{
  \begin{array}{ccc}
    \geq 4 & \textrm{if} & m= 1 \\
    \geq 3m & \textrm{if}& m\geq 2
  \end{array}
\right.$.
\end{thm}

\begin{proof}
Obviously, a plane section of any hypersurface of $\mathbf{A}^N$ with $N\geq 3$ contains at most $\sharp \mathbf{A}^2(\F_2)=4$ points and at most $2$ of them are collinear.
Therefore, since we proved that the $3$ smallest kinds of configurations of points giving low weight codewords are plane configurations, Theorem \ref{main} yields the result.
\end{proof}

%-----------------------------------------------------------------------------
% \noindent \hrulefill Peut \^etre coup\'e \hrulefill
% \subsection{Weight hierarchy of Reed-Muller codes}
% Recall that, from \cite{delsarte} XXXX, the dual of the $q$-ary Reed-Muller code $RM_q(N,r)$ is $RM_q(N,N(q-1)-1-r)$.

% \begin{thm}
%   The reached weights lower than or equal to $3m$ of the Reed--Muller code $RM_q(N,N(q-1)-1-r)$ are: 
% $$\left\{
% \begin{array}{ll}
%   m+2,\ m+3, \ldots, q & \textrm{(only if}\ m+2\leq q\textrm{)} \\
%   2m+2,\ 2m+3, \ldots, 2q & \textrm{(only if}\ m+1\leq q\textrm{)} \\
%   3m & \textrm{(only if}\ m\leq q\textrm{)}
% \end{array}\right. .
% $$
% \end{thm}

% \begin{proof}
%   blablabla
% \end{proof}

% \noindent \hrulefill Fin de coupe \hrulefill
%-----------------------------------------------------------------------------

\section{Key tools for the proof of Theorem \ref{main}}\label{secmlink}

In this section, we state two fundamental results for the proof of Theorem \ref{main} (Lemma \ref{indmin} and Theorem \ref{intersec}).

\subsection{Context} In the present section, $m$ denotes an integer greater than or equal to $1$. 
The base field $k$ is algebraically closed, since Lemma \ref{geom} asserts that treating this case is sufficient.

\subsection{The statements} The following lemma is elementary but very useful in Sections \ref{secsecond} to \ref{secfourth}.

\begin{lem}\label{indmin}
Let $P_1, \ldots, P_s$ be a minimally $m$--linked configuration of points in $\P^r$. Let $d$ and $l$ be two integers satisfying respectively $1 < d < m$ and $1 < l < s$. 
Let $H$ be a hypersurface of degree $d$ containing exactly $l$ of the $P_i$'s. Then, the $s-l$ remaining points are $(m-d)$--linked.
\end{lem}

\begin{proof}
  After a suitable reordering, we have $P_1, \ldots , P_l \in H$ and $P_{l+1}, \ldots ,$ $P_s \notin H$. Assume that $P_{l+1}, \ldots , P_{s}$ are in $(m-d)$--general position. Then, there exists a hypersurface $H'$ of degree $m-d$ containing $P_{l+1}, \ldots , P_{s-1}$ and avoiding $P_s$.
The hypersurface $H \cup H'$ of degree $m$ contains $P_1, \ldots , P_{s-1}$ and avoids $P_s$, which leads to a contradiction thanks to Remark \ref{remmin}. 
\end{proof}

The following statement gives plenty of examples of $m$--linked configurations of points.

\begin{thm}[Cayley--Bacharach]\label{intersec}
  Let $a$ be a positive integer such that $a<m+3$. A family of $a(m+3-a)$ distinct points in $\P^2_k$ lying at the intersection of a curve $C_1$ of degree $a$ and a curve $C_2$ of degree $m+3-a$ having no common irreducible component is minimally $m$--linked.
\end{thm}

\begin{proof}
  Use \cite[Theorem CB4]{EGH} and Remark \ref{remmin0}.
\end{proof}

We conclude the present section by relating the Goppa designed distance for codes on plane curves and Theorem \ref{intersec}.

\subsection{The Goppa designed distance for codes on plane curves}\label{back}
Back to the case of plane curves (see the context in Section \ref{curvecont}).
We proved that the minimum distance of the code $C_L (D,G_m)^{\bot}$ is greater than or equal to the Goppa designed distance which equals $a(m+3-a)$ (see (\ref{deltaG}) page \pageref{deltaG}).
Therefore, the Goppa designed distance for codes on plane curves is closely related to the notion of minimally $m$--linked points in the plane.
In particular, Theorem \ref{intersec} has the following corollary.

\begin{cor}
In the context of Section \ref{curvecont}, assume that the degree $a$ of the plane curve $X$ is greater than or equal to $3$.
If $s=a(m+3-a)$ of the $P_i$'s lie on a curve $Y$ of degree $m+3-a$ which does not contain $X$, then the Goppa designed distance is reached for the code $C_{\Omega}(D,G_m)$.
\end{cor}

\section{First minimal configuration and proof of Theorem \ref{geomain}(\ref{G1})}\label{secfirst}

Obviously, $m+2$ collinear points of a projective space are coplanar and lie at the intersection of a line $L$ and a plane curve $C$ of degree $m+2$ which does not contain $L$.
Therefore, applying Theorem \ref{intersec} for $a=1$, one concludes that $m+2$ collinear points are minimally $m$--linked.
The aim of Proposition \ref{mg1} below is to show that there are no smaller minimally $m$--linked configurations.

\subsection*{Context} In this section the base field $k$ is algebraically closed (it is sufficient to treat this case thanks to Lemma \ref{geom}) and $m\geq 0$ (even if the cases $m=0$ and $1$ are treated in Section \ref{small}, treating them in the present section does not make the proofs longer).

\begin{prop}\label{mg1}
A set of $s \leq m+1$ distinct points $P_1, \ldots , P_s \in \P^r$ is $m$--general.
\end{prop}

\begin{proof}
Let $P_i$ be one of the $s$ points. For all $j\neq i$, there exists a hyperplane  containing $P_j$ and avoiding $P_i$. The union of these $s-1$ hyperplanes is a hypersurface of degree $s-1$ avoiding $P_i$ and containing $P_j$ for all $j\neq i$.
By assumption, $s-1\leq m$. This concludes the proof.
\end{proof}

The following lemma entails the converse statement of Theorem \ref{geomain}(\ref{G1}): \textit{if the minimum distance of $C_L (\Delta, G_m)^{\bot}$ equals $m+2$, then $m+2$ of the $P_i$'s are collinear}.
Moreover, it asserts that the support of a codeword of weight $m+2$ in $C_L (\Delta, G_m)^{\bot}$ is a set of $m+2$ collinear points.

\begin{prop}\label{convmg1}
  Let $P_1, \ldots , P_{m+2}$ be a family of $m$--linked points. Then they are collinear.
\end{prop}

\begin{proof}
Assume that the $P_i$'s are not collinear.
After a suitable reordering of the indexes, $P_m, P_{m+1}$ and $P_{m+2}$ are not collinear and hence there exists a hyperplane $H$ containing $P_{m+1}, P_{m+2}$ and avoiding $P_{m}$.
Therefore, at least $1$ and at most $m$ of the $P_i$'s lie out of $H$ and, from Lemma \ref{indmin}, they are $(m-1)$ linked.
This contradicts Proposition \ref{mg1} applied to $m-1$.
\end{proof}

Let us proceed to the proof of Theorem \ref{geomain}(\ref{G1}).

\begin{proof}[Proof of Theorem \ref{geomain}(\ref{G1})]\label{proofi}
Proposition \ref{mg1} entails that the smallest number of $m$--linked points in a projective space is $\geq m+2$.
Theorem \ref{intersec} entails that $m+2$ collinear points are $m$--linked, which yields the ``if'' part of Theorem \ref{geomain}(\ref{G1}).
The ``only if'' part is a consequence of Proposition \ref{convmg1}.
\end{proof}

\section{Second minimal configuration and proof of Theorem \ref{geomain}(\ref{G2})}\label{secsecond}

\subsection*{Context} In this section, the ambient space is $\P^r$ with $r \geq 2$, the base field $k$ is algebraically closed (see Lemma \ref{geom}) and
$m\geq 2$ (the cases $m=0,1$ have been treated in Section \ref{small}).

\begin{lem}\label{conic}
Let $C$ be a reduced plane conic, $m$ be a positive integer and $P_1, \ldots,$ $P_{2m+2}$ be a family of points of $C$ such that no $m+2$ of them are collinear.
Then, there exists a plane curve $C'$ of degree $m+1$ having no common component with $C$ and intersecting it exactly at the points $P_1, \ldots , P_{2m+2}$. 
\end{lem}

\begin{proof}
  For all $i\in \{1,\ldots, m+1\}$, denote by $L_i$ the line joining $P_i$ and $P_{m+1+i}$. If $C$ is irreducible, then it does not contain any line and the curve $C':=\cup_{i=1}^{m+1} L_i$ is a solution of the problem.
If $C$ is reducible, then, since no $m+2$ of the $P_i$'s are collinear, $C$ is a union of $2$ lines $D_1$ and $D_2$ and each of these lines contains exactly $m+1$ of the $P_i$'s.
After a suitable reordering of the indexes, we have $P_1, \ldots , P_{m+1} \in D_1$ and $P_{m+2}, \ldots , P_{2m+2} \in D_2$. Then, the curve $C':=\cup_{i=1}^{m+1} L_i$ is a solution to the problem.
\end{proof}

From Theorem \ref{intersec} applied to $a=2$ and Lemma \ref{conic}, any set of $2m+2$ points on a plane conic such that no $m+2$ are collinear is minimally $m$--linked.
The purpose of the present section is to prove that there is no other minimally $m$--linked configuration of cardinality $\leq 2m+2$.

\begin{rem}
  It is proved in \cite[Proposition 1]{EGH} that $m+2$ collinear points and $2m+2$ points lying on a conic are the smallest minimally $m$--linked configurations in $\P^2$. However it is not clear that the result holds when the ambient dimension is higher.
\end{rem}

\begin{prop}\label{mg2}
A configuration of $s \leq 2m+1$ distinct points $P_1, \ldots , P_s \in \P^r$ such that no $m+2$ of them are collinear is $m$--general.
\end{prop}

\begin{proof}
For all $m \geq 1$, let $s_m \geq m+2$ be the smallest number of minimally $m$--linked points such that no $m+2$ of them are collinear. From Theorem \ref{intersec} we have $s_m\leq 2m+2$. 
Let us prove that $s_m\geq 2m+2$ by induction on $m$.
\medbreak

\noindent \textbf{Step 1. Initialisation: $\mathbf{m=1}$.}
From Lemma \ref{meq1}, we have $s_1=4$.

\medbreak

\noindent \textbf{Step 2. Induction.} 
Let $m\geq 2$ and assume that $s_{m-1} \geq 2m$.
Let $P_1, \ldots , P_{s_m}$ be a family of minimally $m$--linked points such that no $m+2$ of them are collinear.
Let $c$ be the maximal number of collinear points among $P_1, \ldots , P_{s_m}$. Obviously, we have $2 \leq c$ and, by assumption on the $P_i$'s, we have $c\leq m+1$.

\medbreak

\textbf{Case 2.1. If $\mathbf{c=m+1}$}, then there exists a hyperplane $H$ containing $m+1$ of the $P_i$'s and avoiding all the other ones.
From Lemma \ref{indmin}, the $s_m-m-1$ of the $P_i$'s which lie out of $H$ are $(m-1)$--linked.
Consequently, from Proposition \ref{mg1} we have
$$s_m-m-1 \geq m+1 \quad \textrm{and hence} \quad s_m\geq 2m+2.$$

\medbreak

\textbf{Case 2.2. If $\mathbf{2\leq c \leq m}$}, then, as in the previous step, we prove that $s_m-c$ of the $P_i$'s are $(m-1)$--linked
and, by definition of $c$, no $m+1$ of them are collinear.
By induction, we have 
$$
s_m -c \geq s_{m-1}\geq 2m \quad \textrm{and hence} \quad s_m \geq 2m+2.
$$

\noindent Finally, we always have $s_m \geq 2m+2$.
\end{proof}

Thanks to the previous results we are able to prove a useful and interesting statement asserting that {\it small} minimally $m$--linked configurations are contained in a projective plane.

\begin{prop}\label{coplan1}
For all $m\geq 1$, any minimally $m$--linked configuration of $n\leq 3m$ points is a set of coplanar points.
\end{prop}

\begin{proof}
We prove the result by induction on $m$. If $m=1$, then the result is obvious since any $3$ points are always coplanar.
Let $m> 1$, $n\leq 3m$ and $P_1, \ldots, P_n$ be a minimally $m$--linked configuration of points which we assume to be non--coplanar.
Denote by $s$ the maximal number of coplanar points among them. By assumption, we have $3\leq s<n$.
Moreover, using Proposition \ref{mg2}, one can assume that 
\begin{itemize}
\item[(a)] $n\geq 2m+2$;
\item[(b)] no $m+2$ of the $P_i$'s are collinear.
\end{itemize}

\medbreak

\noindent {\bf Step 1.} Let us prove that $m+1$ of the $P_i$'s are collinear.

\noindent After a suitable reordering, the points $P_1, \ldots , P_s$ are coplanar. Then, there exists a hyperplane $H_0$ containing them and avoiding $P_{s+1}, \ldots , P_n$.
From Lemma \ref{indmin}, the $P_i$'s out of $H_0$ are $(m-1)$--linked.
In particular, $t\leq n-s$ of them are minimally $(m-1)$--linked. After a suitable reordering, $P_{s+1}, \ldots , P_{s+t}$ are minimally $(m-1)$--linked.
Since $s\geq 3$ and thanks to Proposition \ref{mg1}, we get $m+1\leq t \leq 3m-3$. By induction, $P_{s+1}, \ldots , P_{s+t}$ are coplanar.
By definition, $s\geq t \geq m+1$ and $t\leq n-s\leq 2m-1$.
From Proposition \ref{mg2}, the points $P_{s+1}, \ldots, P_{s+t}$ are collinear and $t=m+1$.

\medbreak

\noindent {\bf Step 2.}
Since $m+1$ of the $P_i$'s are collinear and, from (b), no $m+2$ are, there exists a hyperplane $H_1$ containing $m+1$ of them and avoiding all the other ones. From Lemma \ref{indmin}, the points out of $H_1$ are $(m-1)$--linked and their number equals $2m-1$. From the contraposition of Proposition \ref{mg2}, $m+1$ of the points out of $H_1$ are also collinear. Using (a) we split the end of the proof into two cases, both leading to a contradiction.

\smallbreak

{\bf Case 2.1. If} $\mathbf{n>2m+2}$, then there exists a union of two hyperplanes $H_1\cup H_2$ containing $2m+2$ of the $P_i$'s and avoiding the other ones. From Lemma \ref{indmin}, the points out of $H_1 \cup H_2$ are $(m-2)$--linked but their number is $n-(2m+2)\leq m-2$, which contradicts Proposition \ref{mg1} applied to $m-2$.

\smallbreak

{\bf Case 2.2. If} $\mathbf{n=2m+2}$, then the $P_i$'s are contained in a union of two lines $L_1\cup L_2$ which are skew since the $P_i$'s are assumed to be non--collinear. From (b) and after a suitable reordering, $P_1, \ldots P_{m+1}\in L_1$ and $P_{m+2}, \ldots , P_{2m+2}\in L_2$. There exists a hyperplane containing $L_1$ and $P_{2m+2}$. Consequently, from Lemma \ref{indmin}, the points $P_{m+3}, \ldots , P_{2m+2}$ are $(m-1)$--linked contradicting Proposition \ref{mg1} applied to $m-1$. 
\end{proof}

The following proposition yields the converse statement of Theorem \ref{geomain}(\ref{G2}): \textit{if the minimum distance $d$ of $C_L (\Delta, G_m)^{\bot}$ equals $2m+2$, then no $m+2$ of the $P_i$'s are collinear and $2m+2$ of them lie on a plane conic}.
Moreover, the support of a minimum weight codeword of $C_L (\Delta, G_m)^{\bot}$ is contained in a plane conic.

\begin{prop}\label{convmg2}
A minimally $m$--linked configuration of $2m+2$ points such that no $m+2$ of them are collinear is a family of points lying on a plane conic.
\end{prop}

\begin{proof}
From Proposition \ref{coplan1}, the points are coplanar. One concludes using \cite[Proposition 1]{EGH}.
\end{proof}
 
Let us proceed to the proof of Theorem \ref{geomain}(\ref{G2}).

\begin{proof}[Proof of Theorem \ref{geomain}(\ref{G2})]\label{proofii}
From Proposition \ref{mg2}, the smallest number of $m$--linked points such that no $m+2$ are collinear is $\geq 2m+2$.
It is actually an equality from Theorem \ref{intersec} since $2m+2$ points on a plane conic are $m$--linked.
This gives the ``if'' part of the statement.
The ``only if'' part is a consequence of Proposition \ref{convmg2}.
\end{proof}

\section{Third minimal configuration and proof of Theorem \ref{geomain}(\ref{G3})}\label{secthird}

From Theorem \ref{intersec}, we know that $3m$ coplanar points lying at the intersection of a cubic and a curve of degree $m$  having no common component are minimally $m$--linked. The aim of the two remaining sections is to prove that there is no other minimally $m$--linked configuration with cardinality $\leq 3m$.
In this section, we prove that there is no minimally $m$--linked configuration of points of cardinality $<3m$ such that no $m+2$ of the points are collinear and no $2m+2$ of them are on a plane conic.

\subsection*{Context} The ambient space is $\P^r$ with $r\geq 2$, the base field $k$ is algebraically closed and $m\geq 1$ (even if the cases $m=0,1$ have been treated in Section \ref{small}, keeping the case $m=1$ does not make the proofs longer).

First, we need the following elementary lemma.

\begin{lem}\label{cones}
  Let $C$ be a plane conic contained in $\P^r$ and $P_1, \ldots , P_n$ be points avoiding $C$. Then, there exists a hypersurface of degree $2$ containing $C$ and avoiding all the $P_i$'s.
\end{lem}

\begin{proof}
If $r=2$ it is obvious, the expected hypersurface is $C$.
If $r\geq 3$, then consider the set of $3$--codimensional linear subspaces $\Pi \subset \P^r$ such that the cone generated by $C$ over $\Pi$ avoids the $P_i$'s.
One proves easily that this set corresponds to a nonempty open subset of the Grassmanian  $\textrm{Grass}(r-2,k^{r+1})$ (see \cite[Example I.4.1.1]{sch1} for a definition).
\end{proof}

\begin{prop}\label{mg3}
Any $s\leq 3m-1$ distinct points such that no $m+2$ of them are collinear and no $2m+2$ of them lie on a plane conic are $m$--general.
\end{prop}

\begin{proof}
  The method is nearly the same as that of the proof of Proposition \ref{mg2}.
For all $m\geq 1$, denote by $t_m$ the smallest
cardinality of an $m$--linked set of points such that no $2m+2$ of them lie on a plane conic and no $m+2$ of them are collinear.
From Theorem \ref{intersec}, we have $t_m\leq 3m$.
Let us prove that $t_m\geq 3m$ by induction on $m$.

\medbreak

\noindent \textbf{Step 1. Initialisation.}
Proposition \ref{mg2} applied to $m=1$ and $m=2$ respectively entails $t_1 > 3$ and $t_2\geq 6$.
From the same proposition applied to $m=3$,
any $s\leq 7$ points such that no $5$ of them are collinear are $m$--general. Thus, $t_3\geq 8$. 
Moreover, from Proposition \ref{convmg2},
an $8$--tuple of points such that no $5$ of them are collinear and which do not lies on a plane conic is not $m$--linked and hence is $m$--general. Thus, $t_3\geq 9$.

\medbreak

\noindent \textbf{Step 2. Induction.}
Let $m\geq 4$ and $P_1, \ldots , P_{t_m}$ be a minimally $m$--linked configuration of points such that no $m+2$ of them are collinear and no $2m+2$ of them lie on a plane conic, from Proposition \ref{mg2}, we have $t_m\geq 2m+2$.
Moreover, by assumption on the $P_i$'s, since no $2m+2$ of them lie on a plane conic, from Proposition \ref{convmg2}, we have
\begin{equation}\label{am}
t_m > 2m+2
\end{equation}

\noindent Let $c$ be the maximal number of collinear points in $\{P_1, \ldots , P_{t_m}\}$ and $d$ be the maximal number of the $P_i$'s lying on a plane conic.
Obviously, we have
$2 \leq c$.
Moreover, by assumption on the $P_i$'s, we have
$$
c\leq m+1 \quad \textrm{and} \quad d\leq 2m+1.
$$
We consider separately some particular values of $c$ and $d$.

\medbreak

\textbf{Case 2.1. If $\mathbf{d\geq 2m}$,} then let $C$ be a conic containing $d$ of the $P_i$'s.
From Lemma \ref{cones}, there exists a hypersurface of degree $2$ containing $C$ and avoiding all the points out of it.
From Lemma \ref{indmin}, the $t_m-d$ remaining points are $(m-2)$--linked.
Thus, from Proposition \ref{mg1}, we have
$$
t_m-d \geq m \quad \textrm{and hence} \quad t_m \geq 3m. 
$$

\medbreak

\textbf{Case 2.2. If $\mathbf{c = m+1}$ and $\mathbf{d\leq 2m-1}$},
then let $L$ be a line containing $m+1$ of the $P_i$'s.
There exists a hyperplane $H$ containing $L$ and avoiding all $P_i$'s out of $L$.
From Lemma \ref{indmin}, the $t_m-m-1$ of the $P_i$'s lying out of $L$ are $(m-1)$--linked.
Let us consider separately two different situations.
\begin{enumerate}[(a)]
\item If $m+1$ of these points out of $L$ lie on a line $L'$, then one proves by the same reasoning that the $t_m-2m-2$ of the $P_i$'s lying out of $L \cup L'$ are $(m-2)$--linked.
Consequently, from Proposition \ref{mg1} applied to $m-2$, we have
$$
t_m-2m-2 \geq m, \quad \textrm{which entails}\quad t_m \geq 3m+2 \geq 3m.
$$
\item\label{bbb} If no $m+1$ of the points out of $L$ are collinear, then, from Proposition \ref{mg2}, we have
$$
t_m-m-1 \geq 2m \quad \textrm{and hence} \quad t_m \geq 3m+1 \geq 3m.
$$
\end{enumerate}

\medbreak

\textbf{Case 2.3. If $\mathbf{3 \leq c \leq m}$ and $\mathbf{d\leq 2m-1}$,}
then, one proves as in Case 2(\ref{bbb}) that $t_m-c$ of the $P_i$'s are $(m-1)$--linked. 
Moreover, by definition of $c$ and $d$, no $m+1$ of these $t_m-c$ points are collinear and no $2m$ of them lie on a plane conic.
The induction hypothesis yields
$$
t_m-c\geq t_{m-1}\geq 3m-3 \quad \textrm{and hence} \quad t_m \geq 3m.
$$

\medbreak

\textbf{Case 2.4. If $\mathbf{c=2}$, $\mathbf{d\leq 2m-1}$ and the $P_i$'s are not coplanar,} then there exists a hyperplane $H$ containing at least $3$ of the $P_i$'s and avoiding at least $1$ of them. Let $h\geq 3$ be the number of $P_i$'s contained in $H$.
From Lemma \ref{indmin}, the points out of $H$ are $(m-1)$--linked.
Moreover, by assumption, no $m+1$ of them are collinear and no $2m$ of them lie on a plane conic.
By induction, we get
$$
t_m-h \geq t_{m-1} \quad \textrm{and hence}\quad t_m\geq 3m.
$$

\medbreak

From now on, the $P_i$'s are assumed to be coplanar.
Therefore, we always have $d\geq 5$.

\medbreak

\textbf{Case 2.5. If $\mathbf{c=2}$, $\mathbf{2m-2 \leq d \leq 2m-1}$ and the $P_i$'s are coplanar,} then let $C$ be a conic containing $d$ of the $P_i$'s.
From Lemma \ref{indmin}, the  points out of $C$ are $(m-2)$--linked. Since $c=2$ and $m\geq 4$, no $m$ of them are collinear.
Thus, from Proposition \ref{mg2}, we have
$$
t_m - d \geq 2m-2, \quad \textrm{which entails} \quad t_m \geq 4m-4
$$
and, since $m \geq 4$, this entails
$
t_m \geq 3m
$.

\medbreak

\textbf{Case 2.6. If $\mathbf{c=2}$,  $\mathbf{6 \leq d \leq 2m-3}$ and the $P_i$'s are coplanar,} then let $C$ be a conic containing $d$ of the $P_i$'s. From Lemma \ref{indmin}, the points lying out of $C$ are $(m-2)$--linked.
Moreover, by assumption on $c$ and $d$, no $m$ of these points are collinear and no $2m-2$ of them lie on a conic. By the induction hypothesis for $m-2$, we have
$$
t_m-d \geq t_{m-2}\geq 3m-6, \quad \textrm{thus} \quad t_m \geq 3m.
$$

\medbreak

\textbf{Case 2.7. If $\mathbf{c=2}$, $\mathbf{d=5}$ and the $P_i$'s are coplanar but do not lie on a cubic curve,}
then let $C$ be a cubic curve containing at least $9$ of the $P_i$'s.
Such a curve exists since the linear system of plane cubics has dimension $9$.
Denote by $r$ the number of the $P_i$'s contained in $C$. By assumption, $9 \leq r <t_m$ and, from Lemma \ref{indmin}, the $t_m-r$ of the $P_i's$ lying out of $C$ are $(m-3)$--linked.
Moreover, by assumption on $c$ and $d$, no $3$ of these remaining points are collinear and no $6$ of them lie on a cubic. Since $m\geq 4$, we have $(m-3)+2\geq 3$ and $2(m-3)+2 \geq 6$.
Thus, by the induction hypothesis for $m-3$, we have
$$
t_m-r \geq t_{m-3}, \quad \textrm{which entails} \quad t_m\geq 3m-9+r \geq 3m.
$$

\medbreak

\textbf{Case 2.8. If $\mathbf{c=2}$, $\mathbf{d=5}$ and the $P_i$'s lie on a plane cubic curve}, then let $C$ be this cubic curve.
Notice that, by assumption, $m\geq 4$ and, from (\ref{am}), we have $t_m \geq 2m+3$ and hence $t_m \geq 11$.
Since no $3$ of the $P_i$'s are collinear and no $6$ of them lie on a conic and $t_m\geq 11$, one proves easily that $C$ is irreducible.
Then, $t_m \geq 3m$ as a straightforward consequence of Lemma \ref{cubic} below.

\medbreak

\noindent \textbf{Conclusion.} In all the considered cases, we have $t_m \geq 3m$.
\end{proof}

\begin{lem}\label{cubic}
  Let $m$ be an integer greater than or equal to $3$. Let $P_1, \ldots , P_{3m-1}\in \P^2$ be a family of points lying on an irreducible plane cubic curve $C$ such that no $3$ of them are collinear and no $6$ of them lie on a conic. Then, the $P_i$'s are in $m$--general position.
\end{lem}

\begin{proof}
Let $F_C$ be a homogeneous equation of $C$.
Denote by $E_m$ the subspace of $\mathcal{F}_{m,2}$ of homogeneous forms vanishing on $C$ (i.e. $E_m := \mathcal{F}_{m-3,2}F_C$). Choose a subspace $H_m \subset \mathcal{F}_{m,2}$ such that
$
\mathcal{F}_{m,2} = E_m \oplus H_m
$
and let $\Gamma_m$ be the linear system $\P(H_m)$.
It is a linear system of curves of degree $m$ which do \textit{not} contain $C$.
Its dimension is
\begin{equation}\label{dimgamma}
\dim (\Gamma_m)= \dim (\mathcal{F}_{m,2})- \dim (\mathcal{F}_{m-3,2})-1 =3m-1.
\end{equation}

\noindent Let us prove the $m$--generality of $P_1, \ldots , P_{3m-1}$ by induction on $m$.

\medbreak 

\noindent \textbf{Step 1. Initialisation.}
If $m=3$, then consider $8$ points of a plane cubic curve $C$.
Since no $3$ of them are collinear and no $6$ lie on a conic, from \cite[Proposition V.4.3]{H}, the linear system of cubics containing $7$ of them has no other base point. Thus, the points are in $3$-general position.

\medbreak

\noindent \textbf{Step 2. Induction.} Let $m\geq 4$ and assume the induction hypothesis to be true for $m-1$.
 By symmetry on the indexes, to prove the result, it is sufficient to prove the existence of a curve of degree $m$ containing $P_1, \ldots , P_{3m-2}$ and  avoiding $P_{3m-1}$.
We will prove the existence of a curve $D$ of degree $(m-1)$ containing $P_3, \ldots ,$ $P_{3m-2}$, and avoiding $P_{3m-1}$. 
By assumption no $3$ of the $P_i$'s are collinear and hence the line $L$ joining $P_1$ and $P_2$ avoids $P_{3m-1}$.
Consequently, the curve $L \cup D$ of degree $m$ avoids $P_{3m-1}$ and contains all the other $P_i$'s. 

\medbreak

\noindent \textit{Sub-step 2.1} Let $j$ be an integer in $\{1,2,3\}$.
By induction, the points $P_j,P_4, \ldots , P_{3m-2}$ are in $(m-1)$--general position.
Therefore, the maps $\ev_{P_j}, \ev_{P_4}, \ldots , \ev_{P_{3m-2}}$ are linearly independent in $\mathcal{F}_{m-1,2}^{\vee}$ and, since they all vanish on $E_{m-1}$ (recall that $E_{m-1}$ denotes the space of forms of degree $m-1$ vanishing on $C$), they induce independent maps in $(\mathcal{F}_{m-1,2}/E_{m-1})^{\vee} \cong H_{m-1}^{\vee}$.
Let $\Lambda_{j}$ be the maximal sub-system of $\Gamma_{m-1}$ of curves containing $P_j,P_4, \ldots , P_{3m-2}$.
From (\ref{dimgamma}) and since $\ev_{P_j}, \ev_{P_4}, \ldots , \ev_{P_{3m-2}}$ are linearly independent in $H_{m-1}^{\vee}$,we have
$$
\dim (\Lambda_j)= \dim (\Gamma_{m-1})-(3m-4)=0.
$$

\medbreak

\noindent \textit{Sub-step 2.2} For all $j\in \{1,2,3\}$, denote by $D_j$ the single element of $\Lambda_j$.
It is the only element in $\Gamma_{m-1}$ containing the points $P_j, P_4, \ldots , P_{3m-2}$.
For the very same reason, there exists a unique element $D_{3m-1}\in \Gamma_{m-1}$ containing the points $P_{4}, \ldots , P_{3m-1}$.

Let us prove that at least one of the curves $D_1, D_2, D_3$ avoids $P_{3m-1}$.
Assume the negation of the statement, i.e. ``$P_{3m-1}$ lies on $D_1, D_2$ and $D_3$''. Since $D_{3m-1}$ is the unique element of $\Gamma_{m-1}$ containing $P_4, \ldots , P_{3m-1}$, this entails $D_1=D_2=D_3=D_{3m-1}$ and this curve of degree $m-1$ does not contain $C$ and meets it at least at $3m-1$ points.  
But such a situation contradicts B\'ezout's theorem.
Thus, for a suitable ordering of the indexes ${1,2,3}$, the curve $D_3$ avoids $P_{3m-1}$, which concludes the proof.
\end{proof}

% To conclude the present section, let us state the following corollary of Proposition \ref{mg3}, which can be regarded as a generalisation of Corollary \ref{exact}.

% \begin{cor}\label{cardinal}
% Let $m,n$ be two integers satisfying $m \geq 3$ and $n \leq 3m-1$.
% Let $P_1,\ldots , P_n$ be a minimally $m$--linked configuration of points in a projective space, then we have
% either $n=m+2$ and the $P_i$'s are collinear or $n=2m+2$ and the $P_i$'s lie on a plane conic.
% \end{cor}

% \begin{proof}
% Since the $P_i$'s are not in $m$--general position, from Proposition \ref{mg3}, we have either $m+2$ of them are collinear or $2m+2$ of them lie on a plane conic.
% Since the $P_i$'s are \textit{minimally} $m$--linked, they must satisfy one of these two configurations.
% \end{proof}

%-----------------------------------------------------------------------------
% \noindent \hrulefill Peut être coupé \hrulefill

% \begin{proof}[Proof of Theorem \ref{weights}(\ref{W2})]
%   It is almost the same proof as for (\ref{W1}), using respectively Propositions \ref{convmg2} and \ref{mg3} instead of Proposition \ref{convmg1} and Proposition \ref{mg2}.
% \end{proof}
% \noindent \hrulefill Fin de coupe \hrulefill
%-----------------------------------------------------------------------------

\section{End of the proof of Theorem \ref{geomain}}\label{secfourth}

To conclude the proof, it remains to show that the configurations of coplanar points lying at the intersection of a cubic and a degree $m$ curve are the only minimally $m$--linked configurations of cardinality $3m$.

\subsection*{Context} The ambient space is $\P^r$ with $r\geq 2$, the base field $k$ is algebraically closed and $m\geq 3$ (because of Remark \ref{meq2}).

\begin{prop}\label{convmg3}
An $m$--linked configuration of $3m$ points such that no $m+2$ of them are collinear and no $2m+2$ of them lie on a plane conic is a family of coplanar points lying at the intersection of a cubic and a curve of degree $m$ having no common component.
\end{prop}

For the proof of Proposition \ref{convmg3}, we need Lemmas \ref{coconic} and \ref{onacubic}.

\begin{lem}\label{coconic}
Let $n$ be an integer greater than or equal to $6$ and $P_1, \ldots , P_n$ be a family of coplanar points which do not lie on a conic.
Then, there exist $6$ of them which are in $2$--general position.
\end{lem}

\begin{proof}
\textbf{Step 1.}
Let us prove that there exist $5$ of the $P_i$'s which are in $2$--general position.
Proposition \ref{mg2} asserts that $5$ coplanar points are $2$--general if no $4$ of them are collinear.
Since the $P_i$'s do not lie on a conic, they are not collinear.
Therefore, one can reorder the indexes such as $P_1, P_2$ and $P_3$ are not collinear.
For all pairs of distinct integers $i,j\leq n$, denote by $L_{i,j}$ the line joining $P_i$ and $P_j$.
Now we have to prove that there exist two of the $P_i$'s with $i>3$ which do not both lie on one of the lines $L_{1,2}, L_{1,3}$ and $L_{2,3}$.
If not, then the points $P_4, \ldots , P_n$ would all lie on one of the lines $L_{1,2}, L_{1,3}$ and $L_{2,3}$, say $L_{1,2}$.
However, this entails that the $P_i$'s would all lie on the conic $L_{1,2}\cup L_{2,3}$, which yields to a contradiction.

\medbreak

\noindent  \textbf{Step 2.} From the previous step, after a suitable reordering of the indexes, the points $P_1, \ldots , P_5$ are in $2$--general position. Since the linear system of conics in $\P^2$ has dimension $5$, there exists a unique conic $C$ containing $P_1, \ldots , P_5$.
By assumption on the $P_i$'s, $C$ avoids at least one of $P_i$'s, say $P_6$ (after a suitable reordering of the indexes).
Thus, the points $P_1, \ldots ,P_6$ do not lie on a conic.
Finally, this proves that any conic containing $5$ points among $P_1, \ldots , P_6$ avoids the $6$--th one and hence that $P_1, \ldots , P_6$ are in $2$--general position.
\end{proof}

\begin{lem}\label{onacubic}
A minimally $m$--linked  family of $3m$ coplanar points such that no $m+2$ of them are collinear and no $2m+2$ of them lie on a conic, lies on a cubic curve.
\end{lem}

\begin{proof}
Let $P_1, \ldots , P_{3m}$ be such a configuration of points.
To prove the result, we have to treat separately the cases $m=3$ and $4$. 

\medbreak 

\noindent \textbf{Step 1. Small values of $\mathbf{m}$.}
 If $m=3$, then it is obvious since $9$ coplanar points always lie on a cubic.

 If $m=4$, then, since the $P_i$'s are not assumed to be collinear, after a suitable reordering of the indexes, $P_1, P_2$ and $P_3$ are not collinear.
 Let $C$ be a cubic curve containing the points $P_4, \ldots , P_{12}$.
If some of the points $P_1, P_2, P_3$ lie out of $C$, then, from Lemma \ref{indmin}, they are $1$--linked and hence collinear, which yields a contradiction.
Thus, all the $P_i$'s lie on $C$.

If $m=5$, then one can assume that the $P_i$'s are not contained in a conic (if they are, then the result is proved since a conic is contained in plenty of cubics). Lemma \ref{coconic} asserts that $6$ of the $P_i$'s, say $P_1, \ldots, P_6$ are in $2$--general position. Let $C$ be a cubic containing $P_7, \ldots , P_{15}$. If $C$ does not contain all the $P_i$'s, then, from Lemma \ref{indmin}, the $P_i$'s out of $C$ are $2$--linked which contradicts the $2$--generality of $P_1, \ldots , P_6$.

\medbreak

\noindent \textbf{Step 2. For $\mathbf{m\geq 6}$.}
Let $c, d$ be respectively the maximal number of collinear points and of points lying on a conic among the $P_i$'s.

\medbreak

{\bf Case 2.1. If} $\mathbf{d\geq 2m-3}$, then let $Q$ be a conic containing $d$ of the $P_i$'s. From Lemma \ref{indmin}, the $P_i$'s out of $Q$ are $(m-2)$--linked and their number is at most $m+3$. Since $m+3<2m-2$, Proposition \ref{mg2} entails that $m$ of the $P_i$'s out of $Q$ are contained in a line $L$. If $Q\cup L$ contains all the $P_i$'s, then the result is proved. Else, the $P_i$'s out of $Q\cup L$ are $(m-3)$--linked and their number is at most $3$, which contradicts Proposition \ref{mg1}.

\medbreak

{\bf Case 2.2. If} $\mathbf{d = 2m-4}$, then let $Q$ be a conic as in the previous case. The $P_i$'s out of $Q$ are $(m-2)$--linked and their number is $m+4$. If $m\geq 7$, then $2m-2> m+4$ and the result can be obtained by the same manner as in the previous case. If $m=6$, then the $10$ points out of $Q$ cannot lie on a conic since their number is larger than $d=8$. Thus, the $P_i$'s out of $Q$ do not lie on a conic and Proposition \ref{mg3} entails that $m=6$ of the $P_i$'s out of $Q$ are collinear. One can then conclude as in the previous case.

\medbreak

{\bf Case 2.3. If $\mathbf{d<2m-4}$ and $\mathbf{c\geq m-1}$}, then, let $L$ be a line containing at least $m-1$ of the $P_i$'s. From Lemma \ref{indmin}, the $P_i$'s out of $L$ are $(m-1)$--linked and, by assumption on $d$ together with Proposition \ref{mg3}, at least $m+1$ of the $P_i$'s lying out of $L$ are on a line $L'$.
The conic $L\cup L'$ contains at least $2m$ of the $P_i$'s, which contradicts the assumption on $d$.

\medbreak

{\bf Case 2.4. Assume that $\mathbf{d<2m-4}$ and $\mathbf{c<m-1}$}. Let $r$ be the maximal number of the $P_i$'s contained in a cubic. If $r=3m$, then the result is proved. 
Now, assume that $r<3m$. Since the linear system of plane cubics has dimension $9$, we clearly have $r\geq 9$.  Let $C$ be a cubic containing $r$ of the $P_i$'s. From Lemma \ref{indmin}, the $P_i$'s out of $C$ are $(m-3)$--linked.
If $r>9$, then the number of $P_i$'s out of $C$ is $3m-r<3(m-3)$ and, using the assumptions on $c$ and $d$ together with Proposition~\ref{mg3}, these points are in $(m-3)$--general position, which yields a contradiction.

Now, assume that $r=9$.
By induction on $m$ and using the assumptions on $c$ and $d$, the $3(m-3)$ points out $C$ are on a cubic. By definition of $r$, it is possible only if $3(m-3)\leq r=9$, that is $m = 6$ (since $m$ is assumed to be $\geq 6$). From Lemma \ref{indmin}, the $9$ points out of $C$ are $3$--linked. Thus, the linear system of cubics containing these $9$ points has dimension $\geq 1$ and hence, there exists a cubic containing these $9$ points together with a $10$--th one. This contradicts the assumption $r= 9$.
\end{proof}

Now, we can prove Proposition \ref{convmg3}.

\begin{proof}[Proof of Proposition \ref{convmg3}]
Let $P_1, \ldots, P_{3m}$ be an $m$--linked configuration of points such that no $m+2$ of them are collinear and no $2m+2$ lie on a plane conic.
From Proposition \ref{mg3}, these points are actually minimally $m$--linked. 
From Proposition \ref{coplan1}, they are coplanar and from Lemma \ref{onacubic}, they lie on a cubic $C$.
It remains to prove that they lie at the intersection of $C$ with a curve of degree $m$ having no common component with $C$.

To prove this, we will use similar objects as in the proof of Lemma \ref{cubic}. Let $F_C$ be a homogeneous equation of $C$. Let $E_m$ be the subspace of $\mathcal{F}_{m,2}$ of homogeneous forms vanishing on $C$ and let $H_m$ be a complement subspace of $E_m$ in $\mathcal{F}_{m,2}$, that is
$
\mathcal{F}_{m,2}=E_m \oplus H_m.
$
Let $\Gamma_m$ be the linear system $\Gamma_m:=\P (H_m)$. It is a linear system of curves of degree $m$ not containing $C$.
From (\ref{dimgamma}) page \pageref{dimgamma}, we have
$$
\dim (\Gamma_m)=3m-1.
$$
Consequently, there exists an element $D$ of $\Gamma_m$ containing the points $P_1, \ldots ,$ $P_{3m-1}$. Moreover, the curve cannot avoid $P_{3m}$ since the $P_i$'s are minimally $m$--linked. 
It remains to prove that $D$ has no common component with $C$. 

\medbreak

\noindent \textbf{If $\mathbf{C}$ is irreducible}, then it is obvious since the elements of $\Gamma_m$ do not contain $C$.

\medbreak

\noindent \textbf{If $\mathbf{C}$ is reducible}, then $C=C_1 \cup C_2$ such that $C_1$ is a line and $C_2$ a conic (possibly reducible).

First, let us prove that $C_1$ and $C_2$ contain respectively $m$ and $2m$ of the $P_i$'s.
By assumption, at most $m+1$ of the $P_i$'s lie on $C_1$ and at most $2m+1$ of them lie on $C_2$.
If $C_1$ contains $m+1$ of the $P_i$'s, then the $P_i$'s out of it are $(m-1)$--linked and their number is $2m-1$.
Proposition \ref{mg2} entails that $m+1$ of these points are contained in a line $L$ and the $P_i$'s out of $C_1 \cup L$ are $(m-2)$--linked and their number is at most $m-2$, which contradicts Proposition \ref{mg1}. Thus, $C_1$ contains at most $m$ of the $P_i$'s.
If $2m+1$ of the $P_i$'s lie on $C_2$, then from Lemma \ref{indmin}, the $P_i$'s out of $C_2$ are $(m-2)$--linked and their number is $m-1$, which contradicts Proposition \ref{mg1}. Thus $C_2$ contains at most $2m$ of the $P_i$'s.

Finally, after a suitable ordering of the indexes, $P_1, \ldots , P_m \in C_1$ and $P_{m+1}, \ldots , P_{3m} \in C_2$.
Moreover, none of the $P_i$'s lies on $C_1 \cap C_2$.
Suppose that $C_1 \subset D$ and $C_2$ has no common component with $D$.
Then $D=C_1 \cup D_1$ where $D_1$ has degree $m-1$.
Since none of the $P_i$'s lies on $C_1 \cap C_2$, the points $P_{m+1}, \ldots , P_{3m}$ lie on $C_2 \cap D'$, but this contradicts B\'ezout's theorem.

Conversely, if $C_2 \subset D$ and $C_1$ is not contained in $D$, then almost the same reasoning leads also to a contradiction.
\end{proof}

We are now able to conclude the proof of Theorem \ref{geomain} by proving items 
(\ref{G3}) and (\ref{G4}).

\begin{proof}[Proof of Theorem \ref{main}(\ref{T3}) and (\ref{T4})]\label{proofiii}
From Proposition \ref{mg3} the smallest number of $m$--linked points such that no $m+2$ are collinear and no $2m+2$ lie on a plane conic is $\geq 3m$.
From Theorem \ref{intersec}, this inequality is actually an equality since $3m$ points lying at the intersection of two coplanar curves of respective degrees $3$ and $m$ are $m$--linked.
This yields the ``if'' part of Theorem \ref{geomain}(\ref{G3}).
The ``only if'' part is a consequence of Proposition \ref{convmg3}.
Item (\ref{G4}) is a straightforward consequence of (\ref{G1}), (\ref{G2}) and (\ref{G3}).
\end{proof}

\section*{Conclusion}
Using the notion of $m$--generality and in particular that of \textit{being minimally $m$--linked}, we obtain some results on the minimum distance of duals of arbitrary--dimensional algebraic--geometric codes.
For plane curves, these results improve in some situations the well-known Goppa bound.
They also give a method to cleverly puncture such a code on a plane curve in order to drastically increase its dual minimum distance.

From a more geometric point of view, we gave the three smallest configurations of minimally $m$--linked points in any projective space.

To improve Theorem \ref{main} it would be interesting to find further items of this hierarchy.
Notice that these first items correspond to configurations of coplanar points.
Nevertheless, the following ones could correspond, for points in $\P^N$, where $N\geq 3$, to non--coplanar configurations of points.

\section*{Acknowledgements}
The author expresses a deep gratitude to Marc Perret and Daniel Augot for their relevant comments about this article.
Computations in Section \ref{secex} have been made thanks to the software \textsc{Magma}.

\appendix
\section{Varieties not containing plane curves of low degree}\label{append}

From Theorem \ref{main}, to get good codes of the form $C_L(\Delta,G)^{\bot}$, it is interesting to look for varieties which do not contain any plane curves of degree $1,2$ and $3$.
The following result makes possible to check whether a ``generic'' hypersurface of $\P^N$ with fixed degree contains any line, plane conic or plane cubic.
The proof is pretty elementary and uses the same tools as that of \cite[Theorem I.6.4.10]{sch1}. We give it because of a lack of references.

\begin{thm}\label{inclusion}
Let $N, d, r$ be integers such that $N \geq 3$, $d\geq 2$ and $d \geq r \geq 1$.
Then, almost all hypersurfaces of degree $d$ in $\P^N$ do not contain any
plane curve of degree $r$ if
$$
\left(
\begin{array}{c}
d+2
\\ 2
\end{array}
\right)
-
\left(
\begin{array}{c}
d-r+2
\\ 2
\end{array}
\right)-x_{r,N}>0,
$$
where
$$
x_{r,N}=\left\{
\begin{array}{ccccc}
2N-2 & \textrm{if} & r & = & 1\\
\left(
  \begin{array}{c}
    r+2 \\ 2
  \end{array}
\right)+3N-7 & \textrm{if} & r & > & 1
\end{array}
\right. .
$$
\end{thm}
 
\begin{rem}\label{degree4}
  The condition of the theorem is sufficient but not necessary.
For instance, for $N=3$ and $r=3$, we get: \textit{almost all surfaces of $\P^3$ of degree $\geq 5$ do not contain any plane cubic}.
Actually, it is also true for surfaces of degree $4$. Indeed, for $N=3$, $d=4$ and $r=1$ the theorem asserts
that generic surfaces of degree $4$ do not contain any line.
Moreover, it is easy to check that a surface of degree $4$ which does not contain any line cannot contain any plane cubic (consider the plane sections of such a surface).
\end{rem}

\begin{proof}[Proof of Theorem \ref{inclusion}]
\textbf{Notations.} In this proof, for all integers $d,N$, we denote by $\Gamma_{d,N}$ the linear system $\P (\mathcal{F}_{d,N})$ of hypersurfaces of degree $d$ in $\P^N$.
Moreover, for all $r\geq 1$, denote by $X_{r,N}$, the variety parameterising the set of the plane curves of degree $r$ contained in $\P^N$ and by $V_{r,d,N}$ the variety defined by 
$$
V_{r,d,N}:=\{(C,H) \in X_{r,N}\times \Gamma_{d,N}\ |\ C\subset H\}.
$$

\medbreak

\noindent \textbf{Step 1. The variety of lines in $\P^N$: the case $\mathbf{r=1}$.}
For all $N\geq 2$, the variety $X_{1,N}$ is isomorphic to the Grassmanian $\textrm{Grass}(2,k^{N+1})$. Thus,
$$
\dim X_{1,N}=2N-2.
$$
(see \cite[Example I.4.1]{sch1}).

\medbreak

\noindent \textbf{Step 2. The variety of planes curves of degree $\mathbf{r \geq 2}$ in $\P^N$.}
For all $N\geq 2$, the variety parameterising the planes contained in $\P^N$ is isomorphic to the Grassmanian $\textrm{Grass}(3,k^{N+1})$.
This variety has dimension $3N-6$.
Then, for all $r\geq 2$, the variety $X_{r,N}$ is a $\Gamma_{r,2}$--bundle over $\textrm{Grass}(3,k^{N+1})$ and hence has dimension
$$
\dim X_{r,N}=\left(\begin{array}{c} r+2 \\ 2  \end{array}\right)+3N-7.
$$

\medbreak

\noindent 
\textbf{Step 3.} Consider the following diagram
$$
\xymatrix{\relax V_{r,d,N} \ar[rd] \ar@/^/[rrd]^{\varphi_1} \ar@/_/[rdd]_{\varphi_2} & & \\
 & X_{r,N}\times \Gamma_{d,N} \ar[r]^{\pi_1} \ar[d]_{\pi_2}& X_{r,N} \\
 & \Gamma_{d,N} &  }
$$
where $\pi_1$ and $\pi_2$ denote the canonical projections.
To prove the theorem, we have to prove that $\varphi_2$ is not dominant. Thus, it is sufficient to prove that $\dim (V_{r,d,N})<\dim \Gamma_{d,N}$.

Let us compute the dimension of $V_{r,d,N}$.
Notice that, for a given plane curve $C$ of degree $r$ in $\P^N$, the set of hypersurfaces of degree $d$ containing $C$ is parametrised by some projective space $\P^l$ whose dimension $l$ does not depend on $C$. Therefore, $V_{r,d,N}$ is a $\P^l$--bundle over $X_{r,N}$.
Since we know the dimension of $X_{r,N}$, we just have to compute the dimension $l$ of the fibre $F_{r,d,N}$ of $\varphi_1$.

Let $C$ be a plane curve of degree $r$ in $\P^N$ and let $\Pi$ be the plane containing it. Consider the map
$$
\nu : \mathcal{F}_{d,N} \twoheadrightarrow \mathcal{F}_{d,2}
$$
which sends a form of degree $d$ to its restriction to $\Pi$.
The set of forms of degree $d$ in $\mathcal{F}_{d,2}$ vanishing on $C$ is isomorphic to $\mathcal{F}_{d-r,2}$. Therefore, the fibre $F_{r,d,N}$ satisfies
$$
F_{r,d,N} \cong \P (\nu^{-1}(\mathcal{F}_{d-r,2})).
$$
The dimension of $F_{r,d,N}$ is 
$$
\dim F_{r,d,N}=
\dim \mathcal{F}_{d-r,2}+\dim \mathcal{F}_{d,N} -\dim \mathcal{F}_{d,2}  -1.
$$
Finally, we have
$$
\dim \Gamma_{d,N}-\dim V_{r,d,N}=
\left(
\begin{array}{c}
d+2
\\ 2
\end{array}
\right)
-
\left(
\begin{array}{c}
d-r+2
\\ 2
\end{array}
\right)-\dim X_{r,N},
$$
which concludes the proof.
\end{proof}

%\bibliographystyle{abbrv}
%\bibliography{biblio}

\end{document}